\documentclass[12pt]{article}
\usepackage{amssymb}
\usepackage{amsmath}
\usepackage{epsfig}
\usepackage{graphicx}
\usepackage[section]{placeins}
\usepackage[font=small]{caption}

\def\be{\begin{equation}}
\def\ee{\end{equation}}
\def\bq{\begin{eqnarray}}
\def\eq{\end{eqnarray}}

\begin{document}
\begin{center}
{\Large \textbf{Lie Algebra and Conservation Laws for the Time-fractional Heat Equation }}
\end{center}


\begin{center}
{\  Amlan K Halder$^{\dagger}$\footnote{Amlan K Halder, amlanhalder1@gmail.com} , C T Duba$^{\star}$ \footnote{Corresponding author: C T Duba, thama.duba@gmail.com}  and
P G L Leach$^{\mathparagraph,\ddag}$\footnote{PGL Leach, leachp@ukzn.ac.za}\\
$^\dagger$Department of Mathematics,
Pondicherry University, Kalapet - 605014, Puducherry, India\\
$^{\star}$ P.O. Box 11648, Marine Parade - 4056, Durban, South Africa\\
$^\mathparagraph$ Institute For Systems Science, Durban University
of Technology, Durban, South Africa\\
$^\ddag$ School of Mathematics, Statistics and Computer Science, University of KwaZulu - Natal, Durban, South Africa\\

\today}
\end{center}

\begin{abstract}
The Lie symmetry method is applied to derive the point symmetries for the $N$-dimensional fractional heat equation. We find that that the numbers of
symmetries and Lie brackets are reduced significantly as compared to the nonfractional order for all the dimensions. In fact for integer
order linear heat equation the number of solution symmetries is equal to the product of the order and space dimension,
whereas for the fractional case, it is half of the product on the order and space dimension.
We have classified the symmetries and discussed the
Lie algebras and conservational laws.  We generalise the number of symmetries to the $n$-dimensional heat equation.
\end{abstract}

\section{Introduction}
Fractional differential equations (FDEs) are generalisations of differential equations of integer order. The most important advantage of FDEs
in many applications is their nonlocal property. It is known that some real physical phenomena are dependent not only upon their current
state, but also upon their historical (nonlocal) property which can be successfully modelled by using FDEs. Fractional partial differential
equations (FPDEs) appear in various research and engineering applications such as physics, biology, rheology, viscoelasticity, control theory,
signal-processing system-identification and electrochemistry.

Several methods have been developed to find solutions of FPDEs. Some of the methods found in the literature include the exponential function method,
the fractional subequation method, the first integral method, the $G'/G$ expansion method and the Lie symmetry method.

Lie symmetry analysis provides a way of dealing with differential equations and systems thereof. The theory of Lie group analysis is used mainly to
construct similarity reductions and invariant solutions and for obtaining conservation laws. Lie symmetries are also used to reduce the order of the
differential equation as well as the number of independent variables. Lie symmetry analysis has been applied extensively in differential
equations of integral order. Recent studies on Lie symmetry analysis now focus on fractional differential equations. The Lie symmetry analysis of
differential equations has been extended to FDEs by Gazizov {\it et al} (2007) who derived a prolongation formula for the
Riemann-Liouville derivative, enabling one to determine its Lie point symmetries. Others, for example, Wang
{\it et al} (2015) conducted a study on Lie symmetries and conservation laws of time-fractional nonlinear dispersive equation. Lashkarian and
Reza Hejazi (2017) showed that, by using Lie analysis, an FPDE can be reduced to an fractional ordinary differential equation (FODE). Ya\c{s}ar
{\it et al} (2016) and Rui and Zhang (2016) have obtained conservation laws for time-fractional partial differential equations. Wang
{\it et al} (2013), Djordjevic and Atanackovic (2008) and Sahadevan and Bakkyaraj (2012) reduced an FPDE to a nonlinear FODE using the
theory of Sophus Lie.

In our present study we present the classification of symmetries of fractional order arising from a study of a time-fractional heat equation
\begin{equation}
 D_{t}^{\alpha}u=\nabla^{2}u,
\label{eq:heateqn}
\end{equation}
where $0<\alpha\leq 1$ and the Laplacian, $\nabla$, may take any number of space dimensions. Moreover we investigate the Lie algebras and derive
conservation laws arising from the study of heat equation of fractional order. Gazizov {\it et al} (2009) considered a one-dimensional nonlinear
time-fractional diffusion equation and derived the Lie point symmetries and the Lie Brackets using both the Riemann-Liouville
and the Caputo derivatives. Lukashchuk (2015) extended their work by constructing conservation laws for the time-fractional diffusion equation.
We extend this work by considering the N-dimensional linear time-fractional diffusion equation. We
classify
the symmetries, obtain the Lie Brackets and derive conservation laws corresponding to each of the Lie point symmetries of the time-fractional
diffusion equations.

This paper is divided as follows: In Section 2 we give the expression for Riemann-Liouville fractional derivative and provide references wherein
some definitions and properties of the Lie group method to analyse fractional partial differential
equations (FPDEs) are given. In Sections 3, 4 and 5, we present the group analysis of one- to four-dimensional integeral and time-fractional heat
equations.
In Section 7 we present conservation laws of both fractional and integer-valued heat equations and in Section 8 we summarize the results and
conclude.

\section{Preliminaries}
The fractional derivative is given in several forms in the literature. Several definitions of the fractional derivative such as
the Riemann-Liouville, the Gr\"{u}wald-Letnikov, the Weyl, the Caputo, the Riesz and the Miller and Ross have been used by different researchers.
In this article we follow the definition of the Riemann-Liouville fractional derivative given by
\[D_{t}^{\alpha}u(t, x)=
\begin{cases}
 \frac{\partial^{m}u}{\partial t^{m}},& \alpha=m\in N;\\
 \frac{1}{\Gamma(m-\alpha)}\frac{\partial^{m}}{\partial t^{m}}\int^{t}_{0}\frac{u(\tau,x)}{(t-\tau)^{\alpha+1-m}}d\tau, & m-1<\alpha<m,
m\in N.
\end{cases} \]

Symmetry analysis of fractional partial differential equations has been summarized in Sahadevan and Bakkyaraj (2012), Ya\c{s}ar {\it et al}
(2016), Wang {\it et al} (2013) and many others. We explain briefly the method here. The reader can also refer to the articles mentioned before and references therein.\\
 Let a fractional pde be of the form
\begin{equation}\label{2a1}
F(t,x_{i},u,u_{x_{i}},u_{x_{i}x_{i}},u_{t}^{\alpha},u_{x_{i}x_{i}x_{i}},.....)=0,
\end{equation}
 with $0<\alpha<1$, $t$ and $x_{i}$ for $i=1,2,...n$,  represent the independent variables, $u$ as the dependent variable, $u_{x_{i}}$, $u_{x_{i}x_{i}}$, $u_{x_{i}x_{i}x_{i}}$ are the first, second and third derivatives of $u$ with respect to $x_{i}$'s. The fractional pde (\ref{2a1}) is invariant under a one parameter point transformation
\begin{eqnarray}\label{2a2}
\tilde{t}&=&t+\epsilon\xi^{0}(t,x_{i},u)+O(\epsilon^{2}),\nonumber\\
\tilde{x_{i}}&=&x_{i}+\epsilon\xi^{i}(t,x_{i},u)+O(\epsilon^{2}), i=1,2,...n\nonumber\\
\tilde{u}&=&u+\epsilon\eta(t,x_{i},u)+O(\epsilon^{2}),\nonumber\\
\end{eqnarray}
where $\epsilon$ represents an infinitesimal parameter, if and only if
\begin{equation}\label{2a3}
(\tilde{u},\tilde{t},\tilde{x_{i}}){\mid_{{\epsilon=0}}}=(u,t,x_{i}), i=1,2,...n.
\end{equation}
The Lie group $G$ of transformations consisting of the infinitesimal transformations are generated by $$\Gamma=\xi^{0}(t,x_{i},u)\partial_{t}+\xi^{i}(t,x_{i},u)\partial_{x_{i}}+\eta(t,x_{i},u)\partial_{u},$$
 where $\xi^{0}$, $\xi^{i}$ and $\eta$ are the generators of the infinitesimal transformation with respect to $t$, $x_{i}$ and
 $u$. This group $G$ is admitted by the fractional pde (\ref{2a1}).\\
\subsection{The conservation laws}
We mention the preliminaries for the conservation laws of the time fractional heat equation. A vector field, $C=(C^{t}, C^{x})$, where
\begin{eqnarray}
C^{t}&=&C^{t}(t,x,u,u_{x},\cdots),\nonumber\\
C^{x}&=&C^{x}(t,x,u,u_{x},\cdots)
\end{eqnarray}
is called a conserved vector if it satisfies the conservation equation,
\begin{equation}
 D_{t}C^{t}+D_{x}C^{x}=0,
\label{eq:conservationeqn}
\end{equation}
on all solutions of (\ref{eq:heateqn}). Equation(\ref{eq:conservationeqn}) is called a conservation law for (\ref{eq:heateqn})
(Lukashchuk,2015). A conserved vector is called a trivial conserved vector for equation (\ref{eq:heateqn}) if its components $C^{t}$ and
$C^{x}$ vanish on the solution of this equation.

For the time-fractional diffusion equation, (\ref{eq:heateqn}), with the Riemann-Liouville fractional derivative\footnote{ the conservation law for (\ref{eq:heateqn})
with respect to the Caputo derivative is $C^{t}=_{0}I_{t}^{n+1-\alpha}D_{t}u$ (Lukashchuk, 2015)}  is given by
$$ C^{t}=D_{t}^{n-1}(_{0}I_{t}^{n-\alpha}u), \quad C^{x}=-u_{x},\quad n=1,2.$$
Lukashchuk (2015) gives these components as
$$C_{i}^{t}=\phi  (_{0}I_{t}^{1-\alpha}(W_{i}))+J(W_{i},\phi_{t}),$$
$$C_{i}^{x}=\phi_{x}W_{i}-\phi W_{ix},$$
for $\alpha\in(0,1)$, where $i$ corresponds to the number of appropriate symmetries of (\ref{eq:heateqn}).

Hence we need to determine the conserved vectors, $(C^t,C^x)$, in each case.\\

The formal Lagrangian can be introduced as $$L=a(t,x,y,z,w,\cdots)(D^{1-\alpha}_{t}u-\nabla^{2}u),$$ where $a(t,x,y,z,w,\cdots)$ is a new
dependent variable (Ibragimov, 2011; Lukashchuk, 2015; Cao \& Lin, 2014). The Action Integral is $$\int^{T}_{0}\int_{\Omega} L(t,x,u,a,u_{xx})dtdx.$$ Then the Euler-Lagrange equation
(Agrawal, 2002; Atanackovi {\it et al}, 2009) for the fractional
variational principle can be written as
$$\frac{\delta}{\delta u}=\frac{\partial}{\partial u} +D(_{t}^{\alpha})^{*}\frac{\partial}{\partial D^{\alpha}_{t}}u-
D_{x}\frac{\partial}{\partial u_{x}}+D^{2}_{x}\frac{\partial}{\partial u_{xx}},$$
where $(D^{\alpha}_{t})^{*}$ is the adjoint operator of the right Riemann-Liouville derivative
of order $\alpha$ (Ananackovi {\it et al}, 2009).
Therefore, the adjoint equation for the linear heat equation is $$(D_{t}^{\alpha})^{*}a-a_{xx}=0.$$ When
the adjoint equation for the fractional-time heat equation is satisfied for all solutions, $u(t,x)$, upon a substitution $a=\phi(t,x,u)$ such that
$\phi(t,x,u)\neq0$, then it is called nonlinearly self-adjoint.\\
Therefore, the conserved components can be mentioned as \cite{Ibragimov01, Luka01}
\begin{eqnarray}
c^{t}&=&\xi^{0}L+\sum_{m=0}^{n-1} (-1)^{m}D_{t}^{1-\alpha}(W)D_{t}^{m}\frac{\partial L}{u_{t}^{\alpha}}
-(-1)^{n}J\left(W, D_{t}^{n}\frac{\partial L}{u_{t}^{\alpha}}\right)\nonumber\\
\end{eqnarray}\\
where $J$ and $W$ can be defined as
\begin{eqnarray}
J(f,g)&=&\frac{1}{\Gamma(n-\alpha)}\int_{0}^{t}\int_{t}^{T}\frac{f(\tau,x)g(\mu,x)}{(\mu-\tau)^{\alpha+1-n}}d\mu d\tau.\nonumber\\
W&=&\eta-\xi^0 u_t-\sum_{i=1}^{n}\xi^{i} u_{x_{i}},\nonumber\\
\end{eqnarray}
where $\tau$ and $\mu$ are independent variables for $f$ and $g$.  They satisfy the property that
$$D_{t}J(f, g)=f_{t}I^{n-\alpha}_{T}g-g_{0}I_{t}^{n-\alpha}f,$$
where $_{t}I^{n-\alpha}_{T}$ is the right Riemann-Liouville derivative of
order $n-\alpha$. $\phi(t,x)\neq0$ is  nonlinear self-adjoint as defined in Lukshchuk (2015).

\section{Lie point symmetries of the one-dimensional heat equation}
\subsection{Lie point symmetries of the heat equation when $\alpha =1$.}
We perform the Lie symmetry analysis of the one-dimensional heat equation,
\begin{equation}
 u_{t}=u_{xx}.
\label{eq:heq1}
\end{equation}
We obtain the differential operators,
\begin{eqnarray}
\Gamma_{1}&=&\partial_{x},\nonumber\\
\Gamma_{2}&=&2t\partial_{x}-ux\partial_{u},\nonumber\\
\Gamma_{3}&=&\partial_{t},\nonumber\\
\Gamma_{4}&=&2t\partial_{t}+x\partial_{x},\nonumber\\
\Gamma_{5}&=&4t^{2}\partial_{t}+4tx\partial_{x}-u(2t+x^{2})\partial_{u},\nonumber\\
\Gamma_{6}&=&u\partial_{u},\nonumber\\
\Gamma_{7}&=&F(x,t)\partial_{u},
\label{eq:symmetries1}
\end{eqnarray}
where $F(x,t)$ satisfies $$F_{t}-F_{xx}=0.$$

The symmetry $\Gamma_{3}$ is called a time-translational symmetry, $\Gamma_{1}$ and $\Gamma_{2}$ are solution symmetries, whereas $\Gamma_{6}$ is
the homogeneity symmetry
and $\Gamma_{7}$ is the infinite-dimensional symmetry.\\
The nonzero Lie Brackets are

\begin{center}
\begin{tabular}{l l l}
$[\Gamma_{3},\Gamma_{2}]_{LB}=\Gamma_{1},$ & $[\Gamma_{3},\Gamma_{4}]_{LB}=2\Gamma_{3},$\\
$[\Gamma_{3},\Gamma_{5}]_{LB}=-2\Gamma_{6}+4\Gamma_{4}$ & $[\Gamma_{1},\Gamma_{2}]_{LB}=-\Gamma_{6},$\\
$[\Gamma_{1},\Gamma_{5}]_{LB}=2\Gamma_{2},$ & $[\Gamma_{1},\Gamma_{4}]_{LB}=\Gamma_{1},$\\
$[\Gamma_{2},\Gamma_{4}]_{LB}=-\Gamma_{2},$ & $[\Gamma_{4},\Gamma_{5}]_{LB}=2\Gamma_{5}.$\\
\end{tabular}
\end{center}
As it is well known the Lie Algebra is
$$2A_{1}\oplus_{s}sl(2,R)\oplus_{s}A_{1}\oplus_{s}\infty A_{1}.$$

\subsection{Lie point symmetries of the one-dimensional time-fractional heat equation}
We consider the time-fractional heat equation in one-dimensional space
\begin{equation}
 D_{t}^{\alpha}u(x, t)=u_{xx},
\label{eq:fheq1}
\end{equation}
 where $0<\alpha <1$. The Lie symmetries thereof are given by
\begin{eqnarray}
\Gamma_{01}&=&\partial_{x},\nonumber\\
\Gamma_{02}&=&2t\partial_{t}+\alpha x\partial_{x},\nonumber\\
\Gamma_{03}&=&u\partial_{u},\nonumber\\
\Gamma_{04}&=&F(x,t)\partial_{u}.
\label{eq:symmf1}
\end{eqnarray}
The fourth symmetry, $\Gamma_{04}$, is called the infinite-dimensional symmetry, where $F$ is a solution of equation (\ref{eq:fheq1}). It
occurs as a result of the linearity of (\ref{eq:fheq1}) in $u$. We note that for this fractional case we do not have a time-translation symmetry,
$\partial_{t}$. Furthermore, we have only one solution symmetry $\Gamma_{01}$ which is a reduction by one from the nonfractional case. By inspection,
we see
that the Lie Algebra is $A_{1}\oplus A_{2}\oplus\infty A_{1}$.
The Lie Bracket relations among the first three symmetries in equation (\ref{eq:fheq1}) are represented in the table below:
\begin{center}
 \begin{tabular}{c|c|c|c}
  $[,]$ & $\Gamma_{01}$ & $\Gamma_{02}$ & $\Gamma_{03}$ \\
\hline
$\Gamma_{01}$ & $0$ & $0$ & $2\alpha \Gamma_{01}$\\
$\Gamma_{02}$ & $0$ & $0$ & $0$\\
$\Gamma_{03}$ & $-2\alpha \Gamma_{01}$ & $0$ & $0$,
 \end{tabular}
\label{table:LBfheq1}
\end{center}
where each entry, $A_{ij}$, constitutes the Lie Bracket $[\Gamma_{i},\Gamma_{j}]_{LB}$ of two symmetries from equation (\ref{eq:symmf1}) for
$1\leq i, j \leq 3$. We note that $\Gamma_{01}$ generates space translations, $\Gamma_{03}$ reflects the homogeneity property of the equation
(\ref{eq:fheq1}) whereas the $\Gamma_{02}$ is the dilation symmetry. Symmetries $\Gamma_{01}$ to $\Gamma_{03}$ generate the Lie algebra
$A_{1}\oplus A_{2}$.
The Lie Bracket relations satisfied by this three-dimensional algebra are represented in Table \ref{table:LBfheq1}. The nonzero Lie
Bracket relation of the Lie algebra is
$$[\Gamma_{01},\Gamma_{03}]_{LB}=2\alpha \Gamma_{01}.$$By the computation of nonzero Lie Brackets we see that the algebra is the same as what we
obtained from our inspection.

\section{Lie Point symmetries of the two-dimensional heat equation}
\subsection{The case $\alpha=1$.}
For the two-dimensional heat equation,
\begin{equation}
 u_{t}=u_{xx}+u_{yy},
\label{eq:heq2}
\end{equation}
we obtain ten Lie point symmetries given by
\begin{eqnarray}
\Gamma_{21}&=&\partial_{x}, \nonumber\\
\Gamma_{22}&=&\partial_{y},\nonumber\\
\Gamma_{23}&=&2t\partial_{y}-uy\partial_{u}, \nonumber\\
\Gamma_{24}&=&2t\partial_{x}-ux\partial_{u},\nonumber\\
\Gamma_{25}&=&y\partial_{x}-x\partial_{y},\nonumber\\
\Gamma_{26}&=&\partial_{t},\nonumber\\
\Gamma_{27}&=&2t\partial_{t}+x\partial_{x}+y\partial_{y},\nonumber\\
\Gamma_{28}&=&4t^{2}\partial_{t}+4xt\partial_{x}+4yt\partial_{y}-u(4t+x^{2}+
y^{2})\partial_{u},\nonumber\\
\Gamma_{29}&=&u\partial_{u},\nonumber\\
\Gamma_{210}&=&F(x,y,t)\partial_{u}.
\label{eq:symmetries2}
\end{eqnarray}

In this group of symmetries, $\Gamma_{26}$, $\Gamma_{27}$ and $\Gamma_{28}$ constitutes an $sl(2,R)$ subalgebra whereas $\Gamma_{25}$
constitutes an $so(2)$
subalgebra. $\Gamma_{29}$ is the
homogeneity and $\Gamma_{21}$, $\Gamma_{22}$, $\Gamma_{23}$ and $\Gamma_{24}$ are solution symmetries.

The nonzero Lie Brackets are
\begin{center}
 \begin{tabular}{l l l}
  $[\Gamma_{21},\Gamma_{28}]_{LB}=2\Gamma_{24}$ & $[\Gamma_{21},\Gamma_{25}]_{LB}=-\Gamma_{22}$ & $[\Gamma_{21},
\Gamma_{27}]_{LB}=\Gamma_{21}$\\
$[\Gamma_{21},\Gamma_{24}]_{LB}=-\Gamma_{29}$ & $[\Gamma_{28},\Gamma_{26}]_{LB}=-4\Gamma_{29}+\Gamma_{27}$ & $[\Gamma_{28},
\Gamma_{25}]_{LB}=-2\Gamma_{28}$\\
$[\Gamma_{28},\Gamma_{27}]_{LB}=-2\Gamma_{23}$ & $[\Gamma_{26},\Gamma_{27}]_{LB}=2\Gamma_{26}$ & $[\Gamma_{26},\Gamma_{24}]_{LB}=2\Gamma_{21}$\\
$[\Gamma_{26},\Gamma_{23}]_{LB}=2\Gamma_{22}$ & $[\Gamma_{25},\Gamma_{24}]_{LB}=\Gamma_{23}$ & $[\Gamma_{25},\Gamma_{23}]_{LB}=-\Gamma_{24}$\\
$[\Gamma_{25},\Gamma_{22}]_{LB}=\Gamma_{21}$ & $[\Gamma_{27},\Gamma_{24}]_{LB}=2\Gamma_{24}$ & $ [\Gamma_{27},\Gamma_{23}]_{LB}=\Gamma_{23}$\\
$[\Gamma_{27},\Gamma_{22}]_{LB}=-\Gamma_{22}$ & $[\Gamma_{23},\Gamma_{22}]_{LB}=\Gamma_{29}$.
 \end{tabular}
\end{center}
By considering all the subalgebras, the Lie Algebra is $$(4A_{1}\oplus_{s}so(2))\oplus_{s}sl(2,R)\oplus_{s}A_{1}\oplus_{s}\infty A_{1}.$$

\subsection{Lie point symmetries of a two-dimensional time-fractional heat equation}

We consider the two-dimensional time-fractional heat equation:
\begin{equation}
 D_{t}^{\alpha}u(x, y, t)=u_{xx}+u_{yy}.
\label{eq:fheq2}
\end{equation}
In this case for $0<\alpha<1$ the number of symmetries is reduced by four to six. These are given by
\begin{eqnarray}
\Gamma_{11}&=&\partial_{x},\nonumber\\
\Gamma_{12}&=&\partial_{y},\nonumber\\
\Gamma_{13}&=&y\partial_{x}-x\partial_{y},\nonumber\\
\Gamma_{14}&=&4t\partial_{t}+2\alpha x\partial_{x}+2\alpha y\partial_{y}+u(3\alpha-2)\partial_{u},\nonumber\\
\Gamma_{15}&=&u\partial_{u},\nonumber\\
\Gamma_{16}&=&F(x,y,t)\partial_{u},
\label{eq:symmf2}
\end{eqnarray}
where $F(x,y,t)$ satisfies equation (\ref{eq:fheq2}). It can be easily deduced by looking at the symmetries that the Lie Algebra is
$(2A_{1}\oplus_{s}so(2))\oplus_{s}2A_{1}\oplus_{s}\infty A_{1}$.
The nonzero Lie Brackets are
\begin{eqnarray}
& &[\Gamma_{11},\Gamma_{14}]_{LB}=2\alpha \Gamma_{11},\nonumber\\
& &[\Gamma_{14},\Gamma_{13}]_{LB}=-4 \Gamma_{13},\nonumber\\
& &[\Gamma_{14},\Gamma_{12}]_{LB}=-2\alpha \Gamma_{12},\nonumber\\
& &[\Gamma_{11},\Gamma_{16}]_{LB}=G(x,y,t)\partial_{u},
\end{eqnarray}
where $G(x,y,t)$ is an arbitrary solution of equation (\ref{eq:fheq2}).

The algebra is $A_{3,6}\oplus so(2)\oplus_{s}A_{1}\oplus_{s}\infty A_{1}$. It is to be noted that the Lie algebra obtained from nonzero Lie Brackets can be
conflated to the algebra which we have mentioned above. Therefore the Lie algebra obtained by the Lie Brackets can be treated as the composition of the
algebra obtained by inspection.

\section{Lie point symmetries of a three-dimensional heat equation}

\subsection{The case $\alpha=1$}

The three-dimensional heat equation is given by
\begin{equation}
 u_{t}=u_{xx}+u_{yy}+u_{zz}.
\label{eq:heq3}
\end{equation}
The Lie point symmetries of (\ref{eq:heq3}) are given by
\begin{eqnarray}
\Gamma_{31}&=&\partial_{x},\nonumber\\
\Gamma_{32}&=&\partial_{y},\nonumber\\
\Gamma_{33}&=&\partial_{z},\nonumber\\
\Gamma_{34}&=&2t\partial_{y}-uy\partial_{u},\nonumber\\
\Gamma_{35}&=&2t\partial_{x}-ux\partial_{u},\nonumber\\
\Gamma_{36}&=&2t\partial_{z}-uz\partial_{u},\nonumber\\
\Gamma_{37}&=&-y\partial_{x}+x\partial_{y},\nonumber\\
\Gamma_{38}&=&x\partial_{z}-z\partial_{x},\nonumber\\
\Gamma_{39}&=&-z\partial_{y}+y\partial_{z},\nonumber\\
\Gamma_{310}&=&\partial_{t},\nonumber\\
\Gamma_{311}&=&2t\partial_{t}+x\partial_{x}+y\partial_{y}+z\partial_{z},\nonumber\\
\Gamma_{312}&=&4t^{2}\partial_{t}+4xt\partial_{x}+4yt\partial_{y}+4zt\partial_{z}-u(6t+x^{2}+y^{2}+z^{2})\partial_{u},\nonumber \\
\Gamma_{313}&=&u\partial_{u},\nonumber\\
\Gamma_{314}&=&F(x,y,z,t)\partial_{u}.
\label{eq:symmetries3}
\end{eqnarray}
The $sl(2,R)$ symmetries, $\Gamma_{310}$, $\Gamma_{311}$ and $\Gamma_{312}$, remain unchanged. $\Gamma_{37}$, $\Gamma_{38}$ and $\Gamma_{39}$
are the rotation symmetries,
and constitute an $so(3)$ subalgebra,
whereas the $\Gamma_{313}$ is the homogeneity symmetry. We have six solution symmetries given by $\Gamma_{31}$, $\Gamma_{32}$,$\Gamma_{33}$,
$\Gamma_{34}$, $\Gamma_{35}$ and
$\Gamma_{36}$.  $\Gamma_{314}$ is the infinite-dimensional symmetry.

The nonzero Lie Brackets are
\begin{center}
 \begin{tabular}{lll}
  $[\Gamma_{31},\Gamma_{312}]_{LB}=2\Gamma_{35}$, & $[\Gamma_{31},\Gamma_{312}]_{LB}=-\Gamma_{32},$ & $[\Gamma_{31},
\Gamma_{311}]_{LB}=\Gamma_{31}$,\\
$[\Gamma_{31},\Gamma_{35}]_{LB}=-\Gamma_{313}$, & $[\Gamma_{31},\Gamma_{38}]_{LB}=\Gamma_{33},$ & $[\Gamma_{312},
\Gamma_{310}]_{LB}=6\Gamma_{313}-4\Gamma_{311},$\\
$[\Gamma_{312},\Gamma_{311}]_{LB}=-2\Gamma_{312},$ & $[\Gamma_{312},\Gamma_{32}]_{LB}=-2\Gamma_{34}$, & $[\Gamma_{312},
\Gamma_{33}]_{LB}=-2\Gamma_{36},$\\
$[\Gamma_{310},\Gamma_{311}]_{LB}=2\Gamma_{310}$ & $[\Gamma_{310},\Gamma_{35}]_{LB}=2\Gamma_{31}$, & $[\Gamma_{310},
\Gamma_{34}]_{LB}=2\Gamma_{32},$\\
$[\Gamma_{310},\Gamma_{36}]_{LB}=2\Gamma_{33}$ & $[\Gamma_{37},\Gamma_{35}]_{LB}=-\Gamma_{34}$ & $[\Gamma_{37},
\Gamma_{34}]_{LB}=-\Gamma_{35}$\\
$[\Gamma_{37},\Gamma_{32}]_{LB}=-\Gamma_{31}$ & $[\Gamma_{37},\Gamma_{38}]_{LB}=-\Gamma_{39}$ & $[\Gamma_{37},\Gamma_{39}]_{LB}=\Gamma_{38}$\\
$[\Gamma_{311},\Gamma_{35}]_{LB}=\Gamma_{35}$ & $[\Gamma_{311},\Gamma_{34}]_{LB}=\Gamma_{34}$ & $[\Gamma_{311},\Gamma_{32}]_{LB}=-\Gamma_{32}$\\
$[\Gamma_{311},\Gamma_{36}]_{LB}=\Gamma_{36}$ & $[\Gamma_{311},\Gamma_{33}]_{LB}=-\Gamma_{33}$ & $[\Gamma_{35},\Gamma_{38}]_{LB}=\Gamma_{36}$\\
$[\Gamma_{34},\Gamma_{32}]_{LB}=2\Gamma_{313}$ & $[\Gamma_{34},\Gamma_{39}]_{LB}=-\Gamma_{36}$ & $[\Gamma_{32},\Gamma_{39}]_{LB}-\Gamma_{33}$\\
$[\Gamma_{36},\Gamma_{38}]_{LB}=-\Gamma_{35}$ & $[\Gamma_{36},\Gamma_{33}]_{LB}=\Gamma_{313}$ & $[\Gamma_{36},\Gamma_{39}]_{LB}=\Gamma_{34}$\\
$[\Gamma_{38},\Gamma_{33}]_{LB}=-\Gamma_{31}$ & $[\Gamma_{38},\Gamma_{39}]_{LB}=-\Gamma_{37}$ & $[\Gamma_{33},\Gamma_{39}]_{LB}=-\Gamma_{32}$.
\end{tabular}
\end{center}
The Lie Algebra is $$(6A_{1}\oplus_{s} so(3))\oplus_{s} sl(2, R)\oplus_{s}A_{1}\oplus_{s}\infty A_{1}.$$
\subsection{The case $0<\alpha<1$}

The three-dimensional time-fractional heat equation,
\begin{equation}
 D_{t}^{\alpha}u(x, y, z, t)=u_{xx}+u_{yy}+u_{zz},
\label{eq:fheq3}
\end{equation}
has the following Lie point symmetries:
\begin{eqnarray}
\Gamma_{41}&=&\partial_{x},\nonumber\\
\Gamma_{42}&=&\partial_{y},\nonumber\\
\Gamma_{43}&=&\partial_{z},\nonumber\\
\Gamma_{44}&=&-y\partial_{x}+x\partial_{y},\nonumber\\
\Gamma_{45}&=&-y\partial_{z}+z\partial_{y},\nonumber\\
\Gamma_{46}&=&z\partial_{x}-x\partial_{z},\nonumber\\
\Gamma_{47}&=&2t\partial_{t}+\alpha x\partial_{x}+\alpha y\partial_{y}+\alpha z\partial_{z}+u(\alpha-1)\partial_{u},\nonumber\\
\Gamma_{48}&=&u\partial_{u},\nonumber\\
\Gamma_{49}&=&F(x,y,z,t)\partial_{u}.
\label{eq:symmf3}
\end{eqnarray}

Here we lost the $t^{2}$ symmetry.
The nonzero Lie Brackets have been reduced significantly in number and are
\begin{center}
\begin{tabular}{l l}
$[\Gamma_{43},\Gamma_{47}]_{LB}=\alpha \Gamma_{43}$, & $[\Gamma_{43},\Gamma_{45}]_{LB}=-2 \Gamma_{44}$,\\
$[\Gamma_{44},\Gamma_{47}]_{LB}=\alpha \Gamma_{44}$, & $[\Gamma_{44},\Gamma_{45}]_{LB}=-\Gamma_{43},$\\
$[\Gamma_{47},\Gamma_{46}]_{LB}=-2 \Gamma_{46},$ & $[\Gamma_{47},\Gamma_{42}]_{LB}=-\alpha \Gamma_{42}$,\\
$[\Gamma_{47},\Gamma_{41}]_{LB}=-\alpha \Gamma_{41}$ & $[\Gamma_{41},\Gamma_{45}]_{LB}=2 \Gamma_{42}.$
\end{tabular}
\end{center}

The algebra from the Lie Brackets is $$(3A_{1}\oplus_{s}so(3))\oplus_{s}2A_{1}\oplus_{s}\infty A_{1}.$$

\section{Lie point symmetries of the heat equation in four dimensions}

\subsection{The case $\alpha=1$}

The four-dimensional heat equation,
\begin{equation}
 u_{t}=u_{xx}+u_{yy}+u_{zz}+u_{ww},
\label{eq:heq4}
\end{equation}
has the following Lie point symmetries:

\begin{center}
 \begin{tabular}{l l}
$\Gamma_{51}=\partial_{x},$ & $\Gamma_{52}=\partial_{y},$\\
$\Gamma_{53}=\partial_{z},$ & $\Gamma_{54}=\partial_{w},$ \\
$\Gamma_{55}=2t\partial_{y}-uy\partial_{u},$ & $\Gamma_{56}=2t\partial_{x}-ux\partial_{u},$\\
$\Gamma_{57}=2t\partial_{z}-uz\partial_{u},$ & $\Gamma_{58}=2t\partial_{w}-uw\partial_{u},$\\
$\Gamma_{59}=-y\partial_{x}+x\partial_{y},$ & $\Gamma_{510}=y\partial_{w}-w\partial_{y},$ \\
$\Gamma_{511}=-z\partial_{y}+y\partial_{z},$ & $\Gamma_{512}=x\partial_{z}-z\partial_{x},$\\
$\Gamma_{513}=-w\partial_{x}+x\partial_{w},$ & $\Gamma_{514}=z\partial_{w}-w\partial_{z},$ \\
$\Gamma_{515}=\partial_{t},$ & $\Gamma_{516}=2t\partial_{t}+x\partial_{x}+y\partial_{y}+z\partial_{z}+w\partial_{w},$ \\
$\Gamma_{517}=4t^{2}\partial_{t}+4xt\partial_{x}+4yt\partial_{y}+4zt\partial_{z}+4wt\partial_{w}-$\\
$u(8t+x^{2}+y^{2}+z^{2}+
w^{2})\partial_{u},$ \\
$\Gamma_{518}=u\partial_{u},$ &
$\Gamma_{519}=F(x,y,z,w,t)\partial_{u}.$
\label{eq:symmetries4}
\end{tabular}
\end{center}

The nonzero Lie Brackets are given by
\begin{center}
\begin{tabular}{lll}
 $[\Gamma_{51},\Gamma_{517}]_{LB}=2\Gamma_{56}$ & $[\Gamma_{51},\Gamma_{517}]_{LB}=\Gamma_{52}$ & $[\Gamma_{51},\Gamma_{512}]_{LB}=\Gamma_{53}$\\
$[\Gamma_{51},\Gamma_{513}]_{LB}=-\Gamma_{54}$ & $[\Gamma_{51},\Gamma_{516}]_{LB}=\Gamma_{51}$ & $[\Gamma_{51},\Gamma_{56}]_{LB}=-\Gamma_{518}$\\
$[\Gamma_{517},\Gamma_{515}]_{LB}=8\Gamma_{518}-\Gamma_{516}$ & $[\Gamma_{517},\Gamma_{516}]_{LB}=-2\Gamma_{517}$ & $[\Gamma_{517},\Gamma_{52}]_{LB}=-2\Gamma_{58}$\\
$[\Gamma_{517},\Gamma_{53}]_{LB}=-2\Gamma_{57}$ & $[\Gamma_{517},\Gamma_{54}]_{LB}=-2\Gamma_{581}$ & $[\Gamma_{517},\Gamma_{516}]_{LB}=2\Gamma_{515}$\\
$[\Gamma_{515},\Gamma_{56}]_{LB}=2\Gamma_{51}$ & $[\Gamma_{515},\Gamma_{55}]_{LB}=2\Gamma_{53}$ & $[\Gamma_{515},\Gamma_{57}]_{LB}=2\Gamma_{53}$\\
$[\Gamma_{515},\Gamma_{58}]_{LB}=2\Gamma_{54}$ & $[\Gamma_{59},\Gamma_{511}]_{LB}=\Gamma_{512}$ & $[\Gamma_{59},\Gamma_{512}]_{LB}=-\Gamma_{511}$\\
$[\Gamma_{59},\Gamma_{510}]_{LB}=\Gamma_{513}$ & $[\Gamma_{59},\Gamma_{513}]_{LB}=-\Gamma_{510}$ & $[\Gamma_{59},\Gamma_{56}]_{LB}=-\Gamma_{55}$\\
$[\Gamma_{59},\Gamma_{55}]_{LB}=\Gamma_{56}$ & $[\Gamma_{59},\Gamma_{52}]_{LB}=-\Gamma_{51}$ & $[\Gamma_{511},\Gamma_{514}]_{LB}=\Gamma_{510}$\\
$[\Gamma_{511},\Gamma_{512}]_{LB}=\Gamma_{59}$ & $[\Gamma_{511},\Gamma_{510}]_{LB}=-\Gamma_{514}$ & $[\Gamma_{511},\Gamma_{55}]_{LB}=-\Gamma_{57}$\\
$[\Gamma_{511},\Gamma_{57}]_{LB}=\Gamma_{55}$ & $[\Gamma_{511},\Gamma_{52}]_{LB}=\Gamma_{53}$ & $[\Gamma_{511},\Gamma_{53}]_{LB}=\Gamma_{52}$\\
$[\Gamma_{514},\Gamma_{512}]_{LB}=-\Gamma_{513}$ & $[\Gamma_{514},\Gamma_{510}]_{LB}=\Gamma_{511}$ & $[\Gamma_{514},\Gamma_{513}]_{LB}=\Gamma_{512}$\\
$[\Gamma_{514},\Gamma_{57}]_{LB}=\Gamma_{58}$ & $[\Gamma_{514},\Gamma_{58}]_{LB}=-\Gamma_{57}$ & $[\Gamma_{514},\Gamma_{53}]_{LB}=\Gamma_{54}$\\
$[\Gamma_{514},\Gamma_{58}]_{LB}=-\Gamma_{53}$ & $[\Gamma_{512}, \Gamma_{513}]_{LB}=-\Gamma_{514}$ & $[\Gamma_{512},\Gamma_{56}]_{LB}=-\Gamma_{57}$\\
$[\Gamma_{512},\Gamma_{57}]_{LB}=\Gamma_{56}$ & $[\Gamma_{512}, \Gamma_{53}]_{LB}=\Gamma_{51}$ & $[\Gamma_{510},\Gamma_{513}]_{LB}=\Gamma_{59}$\\
$[\Gamma_{510},\Gamma_{55}]_{LB}=\Gamma_{58}$ & $[\Gamma_{510}, \Gamma_{58}]_{LB}=-\Gamma_{55}$ & $[\Gamma_{510},X_{52}]_{LB}=\Gamma_{54}$\\
$[\Gamma_{510},\Gamma_{58}]_{LB}=-\Gamma_{52}$ & $[\Gamma_{513},\Gamma_{54}]_{LB}=-\Gamma_{51}$ & $[\Gamma_{513},\Gamma_{58}]_{LB}=-\Gamma_{56}$\\
$[\Gamma_{513},\Gamma_{56}]_{LB}=\Gamma_{58}$ & $[\Gamma_{516},\Gamma_{54}]_{LB}=-\Gamma_{54}$ & $[\Gamma_{516},\Gamma_{53}]_{LB}=-\Gamma_{53}$\\
$[\Gamma_{516},\Gamma_{52}]_{LB}=-\Gamma_{52}$ & $[\Gamma_{516},\Gamma_{58}]_{LB}=\Gamma_{58}$ & $[\Gamma_{516}, \Gamma_{57}]_{LB}=\Gamma_{57}$\\
$[\Gamma_{516},\Gamma_{55}]_{LB}=\Gamma_{55}$ & $[\Gamma_{516},\Gamma_{56}]_{LB}=\Gamma_{56}$ & $[\Gamma_{55},\Gamma_{52}]_{LB}=\Gamma_{518}$\\
$[\Gamma_{57},\Gamma_{53}]_{LB}=\Gamma_{518}$ & $[\Gamma_{581},\Gamma_{54}]_{LB}=\Gamma_{518}$.\\
\end{tabular}
\end{center}

The Lie Algebra is $$(8A_{1}\oplus_{s}so(4))\oplus_{s}sl(2,R)\oplus_{s}A_{1}\oplus_{s}\infty A_{1}.$$

\subsection{Lie point symmetries of the four-dimensional time-fractional heat equation}

The heat equation in four-dimensional space with fractional time-derivative,
\begin{equation}
 u_{t}^{\alpha}=u_{xx}+u_{yy}+u_{zz}+u_{ww},
\label{eq:fheq4}
\end{equation}
has the following Lie point symmetries
\begin{center}
\begin{tabular} {l l}
$\Gamma_{61}=\partial_{x},$ & $\Gamma_{62}=\partial_{y},$\\
$\Gamma_{63}=\partial_{z},$ & $\Gamma_{64}=\partial_{w},$\\
$\Gamma_{65}=-y\partial_{x}+x\partial_{y},$ & $\Gamma_{66}=-y\partial_{z}+z\partial_{y},$\\
$\Gamma_{67}=y\partial_{w}-w\partial_{y},$ & $\Gamma_{68}=z\partial_{x}-x\partial_{z},$\\
$\Gamma_{69}=-w\partial_{x}+x\partial_{w},$ & $\Gamma_{610}=-w\partial_{z}+z\partial_{w},$\\
$\Gamma_{611}=2t\partial_{t}+\alpha x\partial_{x}+\alpha y\partial_{y}+\alpha z \partial_{z}+\alpha w \partial_{w}+ u (\alpha-1)\partial_{u},$\\
$\Gamma_{612}=u\partial_{u},$ & $\Gamma_{613}=F(x,y,z,w,t)\partial_{u}.$
\end{tabular}
\end{center}
 $\Gamma_{610}$ and $\Gamma_{611}$ are remnants of $sl(2,R)$, $\Gamma_{612}$ is
the homogeneity symmetry,
$\Gamma_{67}$, $\Gamma_{68}$ and $\Gamma_{69}$ constitute an $so(3)$ subalgebra. $\Gamma_{64}$, $\Gamma_{65}$ and $\Gamma_{66}$ are remnants
of $so(4)$ and a representation
of $so(3)$, whereas
$ \Gamma_{61}$, $\Gamma_{62}$ and $\Gamma_{63}$ are solution symmetries. In fact, $\Gamma_{64}$ to $\Gamma_{69}$ are remnants of $so(4)$.

The nonzero Lie Brackets, which are reduced significantly in number compared to the case $\alpha=1$, are
\begin{center}
 \begin{tabular}{lll}
  $[\Gamma_{64},\Gamma_{67}]_{LB}=-\Gamma_{65}$ & $[\Gamma_{64},\Gamma_{68}]_{LB}=-\Gamma_{66}$ & $[\Gamma_{65},\Gamma_{67}]_{LB}=\Gamma_{64}$\\
$[\Gamma_{65},\Gamma_{69}]_{LB}=-\Gamma_{66}$ & $[\Gamma_{62},\Gamma_{67}]_{LB}=\Gamma_{61}$ & $[\Gamma_{62},\Gamma_{69}]_{LB}=\Gamma_{63}$\\
$[\Gamma_{61},\Gamma_{67}]_{LB}=-\Gamma_{62}$ & $[\Gamma_{61},\Gamma_{63}]_{LB}=\Gamma_{68}$ & $[\Gamma_{67},\Gamma_{68}]_{LB}=2\Gamma_{69}$\\
$[\Gamma_{67},\Gamma_{69}]_{LB}=-\Gamma_{68}$ & $[\Gamma_{63},\Gamma_{68}]_{LB}=-\Gamma_{61}$ & $[\Gamma_{63},\Gamma_{69}]_{LB}=-\Gamma_{62}$\\
$[\Gamma_{68},\Gamma_{66}]_{LB}=-\Gamma_{64}$ & $[\Gamma_{68},\Gamma_{69}]_{LB}=\Gamma_{67}$ & $[\Gamma_{66},\Gamma_{69}]_{LB}=\Gamma_{65}$.
 \end{tabular}
\end{center}

The algebra obtained by inspection of the symmetries is $$(4A_{1}\oplus_{s}so(4))\oplus_{s}2A_{1}\oplus_{s}\infty A_{1}.$$

\section{Lie point symmetries of the n-dimensional time-fractional heat equation}

The heat equation in n-dimensional space with fractional time-derivative can be defined as follows,
\begin{equation}
 u_{t}^{\alpha}= u_{x_{1}x_{1}}+u_{x_{2}x_{2}}+u_{x_{3}x_{3}}+..........+u_{x_{n}x_{n}},
\label{eq:fheq4}
\end{equation}
where $x_{1},x_{2},x_{3},.....x_{n}$ are independent variables.
It has the following Lie point symmetries

\begin{eqnarray}
\Gamma_{71}&=&\partial_{x_{i}}, 1\geq i\leq n\nonumber\\
\Gamma_{72}&=&-x_{j}\partial_{x_{i}}+x_{j}\partial_{x_{i}}, i < j\nonumber\\
\Gamma_{73}&=&t\partial_{t}+\alpha \sum_{i=1}^{n} x_{i}\partial_{x_{i}}+ u (\alpha-1)\partial_{u},\nonumber\\
\Gamma_{74}&=&u\partial_{u},\nonumber\\
\Gamma_{75}&=&F(t,x_{1},x_{2},x_{3},.....x_{n})\partial_{u}.\nonumber\\
\end{eqnarray}

\section{The conservation laws for the nonfractional heat equation}

\subsection{Introduction}

We now construct the conservation laws of the heat equation (\ref{eq:heateqn}). The time-fractional diffusion equation (\ref{eq:heateqn}), with the
Riemann-Liouville fractional derivative can be rewritten in the form of conservation laws.

\subsection{Conservation laws for the one-dimensional heat equation}

The Lie point symmetries for the one-dimensional heat equation are given in (\ref{eq:symmetries1}).
The formal Lagrangian can be introduced as $$L=\phi(t,x)(u_{t}-u_{xx}),$$ where $\phi(t,x)$ is a new dependent variable.
The components of the conserved vector for $\Gamma_{1}$ are
\begin{eqnarray}
C^{t}&=&-u_{x}\phi(t,x),\nonumber\\
C^{x}&=&\phi(t,x)(u_{t}-u_{xx})-u_{x}\phi_x+\phi(t,x)u_{xx}.
\end{eqnarray}

For $\Gamma_{2}$ the components of the conserved vector are
\begin{eqnarray}
C^{t}&=&W\phi(t,x),\nonumber\\
C^{x}&=&2t(\phi(t,x)(u_{t}-u_{xx})+W\phi_x+\phi(t,x)(u+xu_{x}+2tu_{xx}),\nonumber\\
W&=&-ux-2tu_{x}.\nonumber\\
\end{eqnarray}

The components of the conserved vectors for $\Gamma_{3}$ are
\begin{eqnarray}
C^{t}&=&\phi(t,x)(u_{t}-u_{xx})+W\phi(t,x),\nonumber\\
C^{x}&=&W\phi_{x}+\phi(t,x)u_{tx},\nonumber\\
W&=&-u_{t}.
\end{eqnarray}

The conserved vectors for $\Gamma_{4}$ are
\begin{eqnarray}
C^{t}&=&2t\phi(t,x)(u_{t}-u_{xx})+W\phi(t,x),\nonumber\\
C^{x}&=&x\phi(t,x)(u_{t}-u_{xx})+W\phi_x+(u_x+xu_{xx}+2tu_{tx}),\nonumber\\
W&=&-2tu_t-xu_x.
\end{eqnarray}

The components of the conserved vectors for $\Gamma_{5}$ are
\begin{eqnarray}
C^{t}&=&4t^2\phi(t,x)(u_{t}-u_{xx})+W\phi(t,x),\nonumber\\
C^{x}&=&4tx\phi(t,x)(u_{t}-u_{xx})+W\phi_x+\phi(t,x)(2tu_{x}+2ux+x^2u_{x}+4t^2u_{tx}+4t(u_{x}+xu_{xx}),\nonumber\\
W&=&-u(2t+x^2)-4t^2u_{t}-4txu_{x}.
\end{eqnarray}

The components of the conserved vectors for $\Gamma_{6}$ are
\begin{eqnarray}
C^{t}&=&W\phi(t,x),\nonumber\\
C^{x}&=&W\phi_x-\phi(t,x)u_{x},\nonumber\\
W&=&u.
\end{eqnarray}
The components of the conserved vectors for $\Gamma_{7}$ are
\begin{eqnarray}
C^{t}&=&W\phi(t,x),\nonumber\\
C^{x}&=&W\phi_x-\phi(t,x)F_{x},\nonumber\\
W&=&F(t,x).\nonumber\\
\end{eqnarray}

\subsection{Conservation laws for the two-dimensional nonfractional heat equation}

The symmetries for the two-dimensional nonfractional heat equation are given in (\ref{eq:symmetries2}).
The formal Lagrangian can be introduced as $$L=\phi(t,x,y)(u_{t}-u_{xx}-u_{yy}),$$ where $\phi(t,x,y)$ is a new dependent
variable.
The components of the conserved vectors for $\Gamma_{21}$ are
\begin{eqnarray}
C^{t}&=&-u_{x}\phi(t,x,y),\nonumber\\
C^{x}&=&\phi(t,x,y)(u_{t}-u_{xx}-u_{yy})+W\phi_{x}+u_{xx},\nonumber\\
C^{y}&=&W\phi_y,\nonumber\\
W&=&-u_x.
\end{eqnarray}

The components of the conserved vectors for $\Gamma_{22}$ are
\begin{eqnarray}
C^{t}&=&-u_y\phi(t,x,y),\nonumber\\
C^{x}&=&W\phi_y,\nonumber\\
C^{y}&=&\phi(t,x,y)(u_{t}-u_{xx}-u_{yy})+W\phi_y+u_{yy},\nonumber\\
W&=&-u_y.
\end{eqnarray}

The components of the conserved vectors for $\Gamma_{23}$ are
\begin{eqnarray}
C^{t}&=&W\phi(t,x,y),\nonumber\\
C^{x}&=&W\phi_{x}+\phi(t,x,y)(y u_{x}),\nonumber\\
C^{y}&=&2t\phi(t,x,y)(u_{t}-u_{xx}-u_{yy})+W\phi_{y}+\phi(u_{y}y+u+2tu_{yy}),\nonumber\\
W&=&-uy-2tu_{y}.
\end{eqnarray}

The components of the conserved vectors for $\Gamma_{24}$ are
\begin{eqnarray}
C^{t}&=&W\phi(t,x,y),\nonumber\\
C^{x}&=&2t\phi(t,x,y)(u_{t}-u_{xx}-u_{yy})+W\phi_{x}+\phi(u+xu_{x}+2tu_{xx}),\nonumber\\
C^{y}&=&W\phi_{y}+\phi(t,x,y)xu_{y},\nonumber\\
W&=&-ux-2tu_{x}.
\end{eqnarray}

The components of the conserved vectors for $\Gamma_{25}$ are
\begin{eqnarray}
C^{t}&=& W\phi(t,x,y),\nonumber\\
C^{x}&=&y\phi(t,x,y)(u_{t}-u_{xx}-u_{yy})+W\phi_{x}-\phi(u_{y}-yu_{xx}),\nonumber\\
C^{y}&=&-x\phi(t,x,y)(u_{t}-u_{xx}-u_{yy})+W\phi_{y}-\phi(xu_{yy}-u_{x}),\nonumber\\
W&=&-yu_{x}+xu_{y}.
\end{eqnarray}

The components of the conserved vectors for $\Gamma_{26}$ are
\begin{eqnarray}
C^{t}&=&\phi(t,x,y)(u_{t}-u_{xx}-u_{yy})+W\phi(t,x,y),\nonumber\\
C^{x}&=&W\phi_x+\phi(t,x,y)u_{tx},\nonumber\\
C^{y}&=&W\phi_y+\phi(t,x,y)u_{ty},\nonumber\\
W&=&-u_{t}.
\end{eqnarray}

The components of the conserved vectors for $\Gamma_{27}$ are
\begin{eqnarray}
C^{t}&=&2t\phi(t,x,y)(u_{t}-u_{xx}-u_{yy})+W\phi(t,x,y),\nonumber\\
C^{x}&=&x\phi(t,x,y)(u_{t}-u_{xx}-u_{yy})+W\phi_x+\phi(2tu_{tx}+xu_{xx}+u_{x}),\nonumber\\
C^{y}&=&y\phi(t,x,y)(u_{t}-u_{xx}-u_{yy})+W\phi_y+\phi(2tu_{ty}+yu_{yy}+u_{y}),\nonumber\\
W&=&-2tu_{t}-xu_{x}-yu_{y}.
\end{eqnarray}

The components of the conserved vectors for $\Gamma_{28}$ are
\begin{align}
C^{t}&=\begin{aligned}[t]
&4t^2\phi(t,x,y)(u_{t}-u_{xx}-u_{yy})+W\phi(t,x,y),\nonumber\\
\end{aligned}\\
C^{x}&=\begin{aligned}[t]
&4xt\phi(t,x,y)(u_{t}-u_{xx}-u_{yy})+W\phi_x+\phi(t,x,y)(4tu_{x}+x^2u_{x}+2xu+y^2u_{x}+2yu+4t^2u_{tx}\nonumber\\
&+4t(u_{x}+xu_{xx})),\nonumber\\
\end{aligned}\\
C^{y}&=\begin{aligned}[t]
&4yt\phi(t,x,y)(u_{t}-u_{xx}-u_{yy})+W\phi_y+\phi(t,x,y)(4tu_{y}+x^2u_{y}+2yu+y^2u_{y}+4t^2u_{ty}\nonumber\\
&+4t(u_{y}+yu_{yy})),\nonumber\\
\end{aligned}\\
W&=\begin{aligned}[t]
&-u(4t+x^2+y^2)-4t^2u_{t}-4xtu_{x}-4ytu_{y}.\nonumber\\
\end{aligned}
\end{align}

The components of the conserved vectors for $\Gamma_{29}$ are
\begin{eqnarray}
C^{t}&=&W\phi(t,x,y),\nonumber\\
C^{x}&=&W\phi_x-\phi(t,x,y)u_{x},\nonumber\\
C^{y}&=&W\phi_y-\phi(t,x,y)u_{y},\nonumber\\
W&=&u.
\end{eqnarray}
The components of the conserved vectors for $\Gamma_{210}$ are
\begin{eqnarray}
C^{t}&=&W\phi(t,x,y),\nonumber\\
C^{x}&=&W\phi_x-\phi(t,x,y)F_{x},\nonumber\\
C^{y}&=&W\phi_y-\phi(t,x,y)F_{y},\nonumber\\
W&=&F(t,x,y).\nonumber\\
\end{eqnarray}

\subsection{Conservation laws of the three-dimensional nonfractional heat equation}

The symmetries for the nonfractional three-dimensional heat equation are given in (\ref{eq:symmetries3}).
The formal Lagrangian can be introduced as $$L=\phi(t,x,y,z)(u_{t}-u_{xx}-u_{yy}-u_{zz}),$$ where $\phi(t,x,y,z)$ is a new dependent
variable.
The components of the conserved vectors for $\Gamma_{31}$ are
\begin{eqnarray}
C^{t}&=&W\phi(t,x,y,z),\nonumber\\
C^{x}&=&\phi(t,x,y,z)(u_{t}-u_{xx}-u_{yy}-u_{zz})+W\phi_x+u_{xx},\nonumber\\
C^{y}&=&W\phi_y,\nonumber\\
C^{z}&=&W\phi_z,\nonumber\\
W&=&-u_{x}.
\end{eqnarray}

The components of the conserved vectors for $\Gamma_{32}$ are
\begin{eqnarray}
C^{t}&=&W\phi(t,x,y,z),\nonumber\\
C^{x}&=&W\phi_x,\nonumber\\
C^{y}&=&\phi(t,x,y,z)(u_{t}-u_{xx}-u_{yy}-u_{zz})+W\phi_y+u_{yy},\nonumber\\
C^{z}&=&W\phi_z,\nonumber\\
W&=&-u_{y}.
\end{eqnarray}

The components of the conserved vectors for $\Gamma_{33}$ are
\begin{eqnarray}
C^{t}&=&W\phi(t,x,y,z),\nonumber\\
C^{x}&=&W\phi_x,\nonumber\\
C^{y}&=&W\phi_y,\nonumber\\
C^{z}&=&\phi(t,x,y,z)(u_{t}-u_{xx}-u_{yy}-u_{zz})+W\phi_z+u_{zz},\nonumber\\
W&=&-u_{z}.
\end{eqnarray}

The components of the conserved vectors for $\Gamma_{34}$ are
\begin{eqnarray}
C^{t}&=&W\phi(t,x,y,z),\nonumber\\
C^{x}&=&W\phi_x+\phi(t,x,y,z)(yu_{y}),\nonumber\\
C^{y}&=&2t\phi(t,x,y,z)(u_{t}-u_{xx}-u_{yy}-u_{zz})+W\phi_y+\phi(u+yu_{y}+2tu_{yy}),\nonumber\\
C^{z}&=&W\phi_z+\phi(yu_{z}),\nonumber\\
W&=&-uy-2tu_{y}.
\end{eqnarray}

The components of the conserved vectors for $\Gamma_{35}$ are
\begin{eqnarray}
C^{t}&=&W\phi(t,x,y,z),\nonumber\\
C^{x}&=&2t\phi(t,x,y,z)(u_{t}-u_{xx}-u_{yy}-u_{zz})+W\phi_x+\phi(u+xu_{x}+2tu_{xx}),\nonumber\\
C^{y}&=&W\phi_y+\phi(t,x,y,z)(xu_y),\nonumber\\
C^{z}&=&W\phi_z+\phi(t,x,y,z)(xu_z),\nonumber\\
W&=&-ux-2tu_{x}.
\end{eqnarray}

The components of the conserved vector for $\Gamma_{36}$ are
\begin{eqnarray}
C^{t}&=&W\phi(t,x,y,z),\nonumber\\
C^{x}&=&W\phi_x+\phi(t,x,y,z)(zu_{x}),\nonumber\\
C^{y}&=&W\phi_y+\phi(t,x,y,z)(zu_{y}),\nonumber\\
C^{z}&=&2t\phi(t,x,y,z)(u_{t}-u_{xx}-u_{yy}-u_{zz})+W\phi_z+\phi(u+zu_{z}+2tu_{zz}),\nonumber\\
W&=&-uz-2tu_{z}.
\end{eqnarray}

The components of the conserved vectors for $\Gamma_{37}$ are
\begin{eqnarray}
C^{t}&=&W\phi(t,x,y,z),\nonumber\\
C^{x}&=&-y\phi(t,x,y,z)(u_{t}-u_{xx}-u_{yy}-u_{zz})+W\phi_x-\phi(yu_{xx}-u_{y}),\nonumber\\
C^{y}&=&x\phi(t,x,y,z)(u_{t}-u_{xx}-u_{yy}-u_{zz})+W\phi_y-\phi(u_{x}-xu_{yy}),\nonumber\\
C^{z}&=&W\phi_z,\nonumber\\
W&=&yu_x-xu_y.
\end{eqnarray}

The components of the conserved vectors for $\Gamma_{38}$ are
\begin{eqnarray}
C^{t}&=&W\phi(t,x,y,z),\nonumber\\
C^{x}&=&-z\phi(t,x,y,z)(u_{t}-u_{xx}-u_{yy}-u_{zz})+W\phi_x-\phi(zu_{xx}-u_{z}),\nonumber\\
C^{y}&=&W\phi_y,\nonumber\\
C^{z}&=&x\phi(t,x,y,z)(u_{t}-u_{xx}-u_{yy}-u_{zz})+W\phi_z-\phi(u_{x}-xu_{zz}),\nonumber\\
W&=&zu_{x}-xu_{z}.
\end{eqnarray}

The components of the conserved vectors for $\Gamma_{39}$ are
\begin{eqnarray}
C^{t}&=&W\phi(t,x,y,z),\nonumber\\
C^{x}&=&W\phi_x,\nonumber\\
C^{y}&=&-z\phi(t,x,y,z)(u_{t}-u_{xx}-u_{yy}-u_{zz})+W\phi_y-\phi(zu_{yy}-u_{z}),\nonumber\\
C^{z}&=&y\phi(t,x,y,z)(u_{t}-u_{xx}-u_{yy}-u_{zz})+W\phi_z-\phi(u_{y}-yu_{zz}),\nonumber\\
W&=&zu_{y}-yu_{z}.
\end{eqnarray}

The components of the conserved vectors for $\Gamma_{310}$ are
\begin{eqnarray}
C^{t}&=&\phi(t,x,y,z)(u_{t}-u_{xx}-u_{yy}-u_{zz})+W\phi(t,x,y,z),\nonumber\\
C^{x}&=&W\phi_x+\phi u_{tx},\nonumber\\
C^{y}&=&W\phi_y+\phi u_{ty},\nonumber\\
C^{z}&=&W\phi_z+\phi u_{tz},\nonumber\\
W&=&-u_{t}.
\end{eqnarray}

The components of the conserved vectors for $\Gamma_{311}$ are
\begin{eqnarray}
C^{t}&=&2t\phi(t,x,y,z)(u_{t}-u_{xx}-u_{yy}-u_{zz})+W\phi(t,x,y,z),\nonumber\\
C^{x}&=&x\phi(t,x,y,z)(u_{t}-u_{xx}-u_{yy}-u_{zz})+W\phi_x+\phi(u_{x}+2tu_{tx}+xu_{xx}),\nonumber\\
C^{y}&=&y\phi(t,x,y,z)(u_{t}-u_{xx}-u_{yy}-u_{zz})+W\phi_y+\phi(u_{y}+2tu_{ty}+yu_{yy}),\nonumber\\
C^{z}&=&z\phi(t,x,y,z)(u_{t}-u_{xx}-u_{yy}-u_{zz})+W\phi_z+\phi(u_{z}+2tu_{tz}+zu_{zz}),\nonumber\\
W&=&-2tu_{t}-xu_{x}-yu_{y}-zu_{z}.
\end{eqnarray}

The components of the conserved vectors for $\Gamma_{312}$ are
\begin{align}
C^{t}&=\begin{aligned}[t]
&4t^2\phi(t,x,y,z)(u_{t}-u_{xx}-u_{yy}-u_{zz})+W\phi(t,x,y,z),\nonumber\\
\end{aligned}\\
C^{x}&=\begin{aligned}[t]
&4xt\phi(t,x,y,z)(u_{t}-u_{xx}-u_{yy}-u_{zz})+W\phi_x+\phi(6tu_x+x^2u_{x}+2xu+y^2u_{x}+z^2u_{x}+4t^{2}u_{tx}\nonumber\\
&+4t(xu_{xx}+u_{x})),\nonumber\\
\end{aligned}\\
C^{y}&=\begin{aligned}[t]
&4yt\phi(t,x,y,z)(u_{t}-u_{xx}-u_{yy}-u_{zz})+W\phi_y+\phi(6tu_{y}+x^2u_{y}+2uy+y^2u_{y}+z^2u_{y}+4t^{2}u_{ty}\nonumber\\
&+4t(u_{y}+yu_{yy})),\nonumber\\
\end{aligned}\\
C^{z}&=\begin{aligned}[t]
&4zt\phi(t,x,y,z)(u_{t}-u_{xx}-u_{yy}-u_{zz})+W\phi_z+\phi(6tu_{z}+x^2u_{z}+y^2u_{z}+z^2u_{z}+2zu+4t^2u_{tz}\nonumber\\
&+4t(u_{z}+zu_{zz})),\nonumber\\
\end{aligned}\\
W&=\begin{aligned}[t]
&-u(6t+x^2+y^2+z^2)-4t^2u_{t}-4xtu_{x}-4ytu_{y}-4ztu_{z}.\nonumber\\
\end{aligned}\\
\end{align}

The components of the conserved vectors for $\Gamma_{313}$ are
\begin{eqnarray}
C^{t}&=&W\phi(t,x,y,z),\nonumber\\
C^{x}&=&W\phi_x-\phi u_{x},\nonumber\\
C^{y}&=&W\phi_y-\phi u_{y},\nonumber\\
C^{z}&=&W\phi_z-\phi u_{z},\nonumber\\
W&=&u.
\end{eqnarray}

The components of the conserved vectors are $\Gamma_{314}$ are
\begin{eqnarray}
C^{t}&=&W\phi(t,x,y,z),\nonumber\\
C^{x}&=&W\phi_x-\phi F_{x},\nonumber\\
C^{y}&=&W\phi_y-\phi F_{y},\nonumber\\
C^{z}&=&W\phi_z-\phi F_{z},\nonumber\\
W&=&F(t,x,y,z).
\end{eqnarray}

\subsection{Conservation laws of the four dimensional non-fractional heat equation}

The formal Lagrangian can be introduced as $$L=\phi(t,x,y,z,w)(u_{t}-u_{xx}-u_{yy}-u_{zz}-u_{ww}),$$ where $\phi(t,x,y,z,w)$ is a new dependent
variable.

The components of the conserved vectors for $\Gamma_{51}$ are
\begin{eqnarray}
C^{t}&=&W\phi(t,x,y,z,w),\nonumber\\
C^{x}&=&\phi(t,x,y,z,w)(u_{t}-u_{xx}-u_{yy}-u_{zz}-u_{ww})+W\phi_x+\phi u_{xx},\nonumber\\
C^{y}&=&W\phi_y,\nonumber\\
C^{z}&=&W\phi_z,\nonumber\\
C^{w}&=&W\phi_w,\nonumber\\
W&=&-u_{x}.
\end{eqnarray}

The components of the conserved vectors for $\Gamma_{52}$ are
\begin{eqnarray}
C^{t}&=&W\phi(t,x,y,z,w),\nonumber\\
C^{x}&=&W\phi_x,\nonumber\\
C^{y}&=&\phi(t,x,y,z,w)(u_{t}-u_{xx}-u_{yy}-u_{zz}-u_{ww})+W\phi_y+\phi u_{yy},\nonumber\\
C^{z}&=&W\phi_z,\nonumber\\
C^{w}&=&W\phi_w,\nonumber\\
W&=&-u_{y}.
\end{eqnarray}

The components of the conserved vectors for $\Gamma_{53}$ are
\begin{eqnarray}
C^{t}&=&W\phi(t,x,y,z,w),\nonumber\\
C^{x}&=&W\phi_x,\nonumber\\
C^{y}&=&W\phi_y,\nonumber\\
C^{z}&=&\phi(t,x,y,z,w)(u_{t}-u_{xx}-u_{yy}-u_{zz}-u_{ww})+W\phi_z+\phi u_{zz},\nonumber\\
C^{w}&=&W\phi_w,\nonumber\\
W&=&-u_{z}.
\end{eqnarray}

The components of the conserved vectors for $\Gamma_{54}$ are
\begin{eqnarray}
C^{t}&=&W\phi_w,\nonumber\\
C^{x}&=&W\phi_x,\nonumber\\
C^{y}&=&W\phi_y,\nonumber\\
C^{z}&=&W\phi_z,\nonumber\\
C^{w}&=&\phi(t,x,y,z,w)(u_{t}-u_{xx}-u_{yy}-u_{zz}-u_{ww})+W\phi_w+\phi u_{ww},\nonumber\\
W&=&-u_{w}.
\end{eqnarray}

The components of the conserved vectors for $\Gamma_{55}$ are
\begin{eqnarray}
C^{t}&=&W\phi(t,x,y,z,w),\nonumber\\
C^{x}&=&W\phi_x+\phi(t,x,y,z,w)(yu_{x}),\nonumber\\
C^{y}&=&2t\phi(t,x,y,z,w)(u_{t}-u_{xx}-u_{yy}-u_{zz})+W\phi_y+\phi(u+yu_{y}+2tu_{yy}),\nonumber\\
C^{z}&=&W\phi_z+\phi(yu_{z}),\nonumber\\
C^{w}&=&W\phi_w+\phi(yu_{w}),\nonumber\\
W&=&-uy-2tu_{y}.
\end{eqnarray}

The components of the conserved vectors for $\Gamma_{56}$ are
\begin{eqnarray}
C^{t}&=&W\phi(t,x,y,z,w),\nonumber\\
C^{x}&=&2t\phi(t,x,y,z,w)(u_{t}-u_{xx}-u_{yy}-u_{zz})+W\phi_x+\phi(u+xu_{x}+2tu_{xx}),\nonumber\\
C^{y}&=&W\phi_y+\phi(t,x,y,z,w)(xu_y),\nonumber\\
C^{z}&=&W\phi_z+\phi(t,x,y,z,w)(xu_z),\nonumber\\
C^{w}&=&W\phi_w+\phi(t,x,y,z,w)(xu_w),\nonumber\\
W&=&-ux-2tu_{x}.
\end{eqnarray}

The components of the conserved vector for $\Gamma_{57}$ are
\begin{eqnarray}
C^{t}&=&W\phi(t,x,y,z,w),\nonumber\\
C^{x}&=&W\phi_x+\phi(t,x,y,z,w)(zu_{x}),\nonumber\\
C^{y}&=&W\phi_y+\phi(t,x,y,z,w)(zu_{y}),\nonumber\\
C^{z}&=&2t\phi(t,x,y,z,w)(u_{t}-u_{xx}-u_{yy}-u_{zz})+W\phi_z+\phi(u+zu_{z}+2tu_{zz}),\nonumber\\
C^{w}&=&W\phi_w+\phi(t,x,y,z,w)(zu_{w}),\nonumber\\
W&=&-uz-2tu_{z}.
\end{eqnarray}

The components of the conserved vector for $\Gamma_{58}$ are
\begin{eqnarray}
C^{t}&=&W\phi(t,x,y,z,w),\nonumber\\
C^{x}&=&W\phi_x+\phi(t,x,y,z,w)(wu_{x}),\nonumber\\
C^{y}&=&W\phi_y+\phi(t,x,y,z,w)(wu_{y}),\nonumber\\
C^{z}&=&W\phi_z+\phi(t,x,y,z,w)(wu_{z}),\nonumber\\
C^{w}&=&2t\phi(t,x,y,z,w)(u_{t}-u_{xx}-u_{yy}-u_{zz})+W\phi_w+\phi(u+wu_{w}+2tu_{ww}),\nonumber\\
W&=&-uw-2tu_{w}.
\end{eqnarray}

The components of the conserved vectors for $\Gamma_{59}$ are
\begin{eqnarray}
C^{t}&=&W\phi(t,x,y,z,w),\nonumber\\
C^{x}&=&-y\phi(t,x,y,z,w)(u_{t}-u_{xx}-u_{yy}-u_{zz})+W\phi_x-\phi(yu_{xx}-u_{y}),\nonumber\\
C^{y}&=&x\phi(t,x,y,z,w)(u_{t}-u_{xx}-u_{yy}-u_{zz})+W\phi_y-\phi(u_{x}-xu_{yy}),\nonumber\\
C^{z}&=&W\phi_z.\nonumber\\
C^{w}&=&W\phi_w,\nonumber\\
W&=&yu_x-xu_y.
\end{eqnarray}
The components of the conserved vectors for $\Gamma_{510}$ are
\begin{eqnarray}
C^{t}&=&W\phi(t,x,y,z,w),\nonumber\\
C^{x}&=&W\phi_x,\nonumber\\
C^{y}&=&-w\phi(t,x,y,z,w)(u_{t}-u_{xx}-u_{yy}-u_{zz}-u_{ww})+W\phi_y-\phi(wu_{yy}-u_{w}),\nonumber\\
C^{z}&=&W\phi_z,\nonumber\\
C^{w}&=&y\phi(t,x,y,z,w)(u_{t}-u_{xx}-u_{yy}-u_{zz}-u_{ww})+W\phi_w-\phi(u_y-yu_{ww}),\nonumber\\
W&=&wu_{y}-yu_{w}.
\end{eqnarray}

The components of the conserved vectors for $\Gamma_{511}$ are
\begin{eqnarray}
C^{t}&=&W\phi(t,x,y,z,w),\nonumber\\
C^{x}&=&-z\phi(t,x,y,z,w)(u_{t}-u_{xx}-u_{yy}-u_{zz})+W\phi_x-\phi(zu_{xx}-u_{z}),\nonumber\\
C^{y}&=&W\phi_y,\nonumber\\
C^{z}&=&x\phi(t,x,y,z,w)(u_{t}-u_{xx}-u_{yy}-u_{zz})+W\phi_z-\phi(u_{x}-xu_{zz}),\nonumber\\
C^{w}&=&W\phi_w,\nonumber\\
W&=&zu_{x}-xu_{z}.
\end{eqnarray}

The components of the conserved vectors for $\Gamma_{512}$ are
\begin{eqnarray}
C^{t}&=&W\phi(t,x,y,z,w),\nonumber\\
C^{x}&=&W\phi_x,\nonumber\\
C^{y}&=&-z\phi(t,x,y,z,w)(u_{t}-u_{xx}-u_{yy}-u_{zz})+W\phi_y-\phi(zu_{yy}-u_{z}),\nonumber\\
C^{z}&=&y\phi(t,x,y,z,w)(u_{t}-u_{xx}-u_{yy}-u_{zz})+W\phi_z-\phi(u_{y}-yu_{zz}),\nonumber\\
C^{w}&=&W\phi_w,\nonumber\\
W&=&zu_{y}-yu_{z}.
\end{eqnarray}

The components of the conserved vector for $\Gamma_{513}$ are
\begin{eqnarray}
C^{t}&=&W\phi(t,x,y,z,w),\nonumber\\
C^{x}&=&-w\phi(t,x,y,z,w)(u_{t}-u_{xx}-u_{yy}-u_{zz}-u_{ww})+W\phi_x-\phi(wu_{xx}-u_{w},\nonumber\\
C^{y}&=&W\phi_y,\nonumber\\
C^{z}&=&W\phi_z,\nonumber\\
C^{w}&=&x\phi(t,x,y,z,w)(u_{t}-u_{xx}-u_{yy}-u_{zz}-u_{ww})+W\phi_w-\phi(u_{x}-xu_{ww}),\nonumber\\
W&=&wu_{x}-xu_{w}.
\end{eqnarray}

The components of the conserved vector for $\Gamma_{514}$ are
\begin{eqnarray}
C^{t}&=&W\phi(t,x,y,z,w),\nonumber\\
C^{x}&=&W\phi_x,\nonumber\\
C^{y}&=&W\phi_y,\nonumber\\
C^{z}&=&-w\phi(t,x,y,z,w)(u_{t}-u_{xx}-u_{yy}-u_{zz}-u_{ww})+W\phi_z-\phi(wu_{zz}-u_{w}),\nonumber\\
C^{w}&=&z\phi(t,x,y,z,w)(u_{t}-u_{xx}-u_{yy}-u_{zz}-u_{ww})+W\phi_w-\phi(u_z-zu_{ww}),\nonumber\\
W&=&wu_{z}-zu_{w}.
\end{eqnarray}

The components of the conserved vector for $\Gamma_{515}$ are
\begin{eqnarray}
C^{t}&=&\phi(t,x,y,z,w)(u_{t}-u_{xx}-u_{yy}-u_{zz}-u_{ww})+W\phi(t,x,y,z,w),\nonumber\\
C^{x}&=&W\phi_x+\phi u_{tx},\nonumber\\
C^{y}&=&W\phi_y+\phi u_{ty},\nonumber\\
C^{z}&=&W\phi_z+\phi u_{tz},\nonumber\\
C^{w}&=&W\phi_w+\phi u_{tw},\nonumber\\
W&=&-u_{t}.
\end{eqnarray}

The components of the conserved vector for $\Gamma_{516}$ are
\begin{eqnarray}
C^{t}&=&2t\phi(t,x,y,z,w)(u_{t}-u_{xx}-u_{yy}-u_{zz}-u_{ww})+W\phi(t,x,y,z,w),\nonumber\\
C^{x}&=&x\phi(t,x,y,z,w)(u_{t}-u_{xx}-u_{yy}-u_{zz}-u_{ww})+W\phi_x+\phi(u_{x}+2tu_{tx}+u_{x}+xu_{xx}),\nonumber\\
C^{y}&=&y\phi(t,x,y,z,w)(u_{t}-u_{xx}-u_{yy}-u_{zz}-u_{ww})+W\phi_y+\phi(u_{y}+2tu_{ty}+u_{y}+yu_{yy}),\nonumber\\
C^{z}&=&z\phi(t,x,y,z,w)(u_{t}-u_{xx}-u_{yy}-u_{zz}-u_{ww})+W\phi_z+\phi(u_{z}+2tu_{tz}+u_{z}+zu_{zz}),\nonumber\\
C^{w}&=&w\phi(t,x,y,z,w)(u_{t}-u_{xx}-u_{yy}-u_{zz}-u_{ww})+W\phi_w+\phi(u_{z}+2tu_{tz}+u_{w}+zu_{zz}),\nonumber\\
W&=&-2tu_{t}-xu_{x}-yu_{y}-zu_{z}-wu_{w}.
\end{eqnarray}

The components of the conserved vectors for $\Gamma_{517}$ are
\begin{align}
C^{t}&=\begin{aligned}[t]
&4t^2\phi(t,x,y,z,w)(u_{t}-u_{xx}-u_{yy}-u_{zz}-u_{ww})+W\phi(t,x,y,z,w),\nonumber\\
\end{aligned}\\
C^{x}&=\begin{aligned}[t]
&4xt\phi(t,x,y,z,w)(u_{t}-u_{xx}-u_{yy}-u_{zz}-u_{ww})+W\phi_x\nonumber\\
&+\phi(8tu_x+x^2u_{x}+2xu+y^2u_{x}+z^2u_{x}+4t^{2}u_{tx}+4t(xu_{xx}+u_{x})),\nonumber\\
\end{aligned}\\
C^{y}&=\begin{aligned}[t]
&4yt\phi(t,x,y,z,w)(u_{t}-u_{xx}-u_{yy}-u_{zz}-u_{ww})+W\phi_y\nonumber\\
&+\phi(8tu_{y}+x^2u_{y}+2uy+y^2u_{y}+z^2u_{y}+4t^{2}u_{ty}+4t(u_{y}+yu_{yy})),\nonumber\\
\end{aligned}\\
C^{z}&=\begin{aligned}[t]
&4zt\phi(t,x,y,z,w)(u_{t}-u_{xx}-u_{yy}-u_{zz}-u_{ww})+W\phi_z\nonumber\\
&+\phi(8tu_{z}+x^2u_{z}+y^2u_{z}+z^2u_{z}+2zu+4t^2u_{tz}+4t(u_{z}+zu_{zz})),\nonumber\\
\end{aligned}\\
C^{w}&=\begin{aligned}[t]
&4wt\phi(t,x,y,z,w)(u_{t}-u_{xx}-u_{yy}-u_{zz}-u_{ww})+W\phi_w\nonumber\\
&+\phi(8tu_{w}+x^2u_{w}+y^2u_{w}+z^2u_{w}+2wu+4t^2u_{tw}+4t(u_{w}+wu_{ww})),\nonumber\\
\end{aligned}\\
W&=\begin{aligned}[t]
&-u(8t+x^2+y^2+z^2+w^2)-4t^2u_{t}-4xtu_{x}-4ytu_{y}-4ztu_{z}-4wtu_{w}.\nonumber\\
\end{aligned}\\
\end{align}

The components of the conserved vectors for $\Gamma_{518}$ are
\begin{eqnarray}
C^{t}&=&W\phi(t,x,y,z,w),\nonumber\\
C^{x}&=&W\phi_x-\phi u_{x},\nonumber\\
C^{y}&=&W\phi_y-\phi u_{y},\nonumber\\
C^{z}&=&W\phi_z-\phi u_{z},\nonumber\\
C^{w}&=&W\phi_w-\phi u_{w},\nonumber\\
W&=&u.
\end{eqnarray}

The components of the conserved vectors are $\Gamma_{519}$ are
\begin{eqnarray}
C^{t}&=&W\phi(t,x,y,z,w),\nonumber\\
C^{x}&=&W\phi_x-\phi F_{x},\nonumber\\
C^{y}&=&W\phi_y-\phi F_{y},\nonumber\\
C^{z}&=&W\phi_z-\phi F_{z},\nonumber\\
C^{w}&=&W\phi_w-\phi F_{w},\nonumber\\
W&=&F(t,x,y,z,w).
\end{eqnarray}

\section{Conservation laws with the Riemann-Liouville fractional derivative for one-dimensional heat equation}

The conservation laws with respect to each $\Gamma_{0i}$, $i=1,2,3,4$ are as follows:
For $\Gamma_{01}$ the conserved vectors are
\begin{eqnarray}
C^t&=&\phi(t,x)D_{t}^{1-\alpha}W+J(W,\phi_{t}),\nonumber\\
C^x&=&\phi(t,x)(D_{t}^{\alpha}u-u_{xx})+W\phi_{x}-\phi(t,x)u_{xx},\nonumber\\
W&=&-u_{x},
\end{eqnarray}

For $\Gamma_{02}$ the conserved vectors are
\begin{eqnarray}
C^t&=&2t\phi(t,x)(D_{t}^{\alpha}u-u_{xx})+\phi(t,x)D_{t}^{1-\alpha}(W)+J(W,\phi_{t}),\nonumber\\
C^x&=&\alpha x \phi(t,x)(D_{t}^{\alpha}u-u_{xx})+W\phi_{x}-\phi(t,x)(2tu_{xt}-\alpha x u_{xx}),\nonumber\\
W&=&2tu_t-\alpha x u_{x}.
\end{eqnarray}
For $\Gamma_{03}$ we have the following components of the conserved vector
\begin{eqnarray}
C^t&=&\phi(t,x)D_{t}^{1-\alpha}W+J(W,\phi_{t}),\nonumber\\
C^x&=&W\phi_{x}-2\phi(t,x)u_{x},\nonumber\\
W&=&u.
\end{eqnarray}
For $\Gamma_{04}$ the conserved vectors are
\begin{eqnarray}
C^t&=&\phi(t,x)D_{t}^{1-\alpha}F(t,x)+J(F(t,x),\phi_{t}),\nonumber\\
C^x&=&F(t,x)\phi_{x}-\phi(t,x)F_{x}.
\end{eqnarray}

\subsection{The conservation laws for the two-dimensional time-fractional heat equation}

The symmetries are given in (\ref{eq:symmf2}).
For $\Gamma_{11}$ the components of the conserved vectors are
\begin{eqnarray}
C^t&=&D_{t}^{1-\alpha}(W)\phi(t,x,y)+J(W,\phi_t),\nonumber\\
C^x&=&\phi(t,x,y)(D_{t}^{\alpha}u-u_{xx}-u_{yy})+W\phi_{x}+u_{xx}\phi,\nonumber\\
C^y&=&W\phi_y+\phi u_{xy},\nonumber\\
W&=&-u_{x}.
\end{eqnarray}
For $\Gamma_{12}$ the components of the conserved vectors are
\begin{eqnarray}
C^t&=&D_{t}^{1-\alpha}(W)\phi(t,x,y)+J(-u_y,\phi_t),\nonumber\\
C^x&=&W\phi_{x}+\phi u_{xy},\nonumber\\
C^y&=&\phi(t,x,y)(D_{t}^{\alpha}u-u_{xx}-u_{yy})+W\phi_y+u_{yy}\phi(t,x,y),\nonumber\\
W&=&-u_{y}.
\end{eqnarray}
For $\Gamma_{13}$ the components of the conserved vectors are
\begin{eqnarray}
C^t&=&D_{t}^{1-\alpha}(W)\phi(t,x,y)+J(W,\phi_t),\nonumber\\
C^x&=& y\phi(t,x,y)(D_{t}^{\alpha}u-u_{xx}-u_{yy})+W\phi_x-\phi(t,x,y)(u_{y}+xu_{xy}-yu_{xx}),\nonumber\\
C^y&=&-x\phi(t,x,y)(D_{t}^{\alpha}u-u_{xx}-u_{yy}) +W\phi_y-\phi(t,x,y)(xu_{yy}-(yu_{xy}+u_{x}),\nonumber\\
W&=&-yu_{x}+xu_{y}.
\end{eqnarray}
For $\Gamma_{14}$ the components are
\begin{eqnarray}
C^t&=&4t(D_{t}^{\alpha}u-u_{xx}-u_{yy})+D_{t}^{1-\alpha}(W)\phi(t,x,y)+J(W,\phi_t),\nonumber\\
C^x&=& 2\alpha x\phi(t,x,y)(D_{t}^{\alpha}u-u_{xx}-u_{yy})+W\phi_x-(3\alpha u_{x}-2u_{x}-4tu_{xt}-2\alpha(u_{x}+x u_{xx})-2\alpha y u_{xy})\phi(t,x,y),\nonumber\\
C^y&=&2\alpha y\phi(t,x,y)(D_{t}^{\alpha}u-u_{xx}-u_{yy})+W\phi_y+(3\alpha u_{y}-2 u_{y}-4t u_{yt}-2\alpha x u_{xy}-2\alpha(y u_{yy}+u_{y}))\phi(t,x,y),\nonumber\\
W&=&u(3\alpha-2)-4tu_t-2\alpha x u_{x}-2\alpha y u_{y}.
\end{eqnarray}
For $\Gamma_{15}$ the components are
\begin{eqnarray}
C^t&=&D_{t}^{1-\alpha}(W)\phi(t,x,y)+J(u,\phi_t),\nonumber\\
C^x&=&u\phi_x-u_x\phi(t,x,y),\nonumber\\
C^y&=&u\phi_y-u_y\phi(t,x,y),\nonumber\\
W&=&u.\nonumber\\
\end{eqnarray}
For $\Gamma_{16}$ the conserved vectors are
\begin{eqnarray}
C^t&=&D_{t}^{1-\alpha}(F(t,x,y))\phi(t,x,y)+J(F(t,x,y),\phi_t),\nonumber\\
C^x&=&F(t,x,y)\phi_x-F_{x}\phi(t,x,y),\nonumber\\
C^y&=&F(t,x,y)\phi_y-F_y\phi(t,x,y).\nonumber\\
\end{eqnarray}

\subsection{The conservation laws for the three-dimensional time-fractional heat equation}

The symmetries are given in (\ref{eq:symmf3}).
The conserved vectors with respect to each symmetry are as follows:
For $\Gamma_{41}$ the components of the conserved vector are
\begin{eqnarray}
C^{t}&=& D_{t}^{1-\alpha}(W)\phi(t,x,y,z)+J(W,\phi_t),\nonumber\\
C^{x}&=&\phi(t,x,y,z)(D_{t}^{1-\alpha}u-(u_{xx}+u_{yy}+u_{zz}))+W\phi_x+\phi(t,x,y,z)u_{xx},\nonumber\\
C^{y}&=&W\phi_y+\phi(t,x,y,z)u_{xy},\nonumber\\
C^{z}&=&W\phi_z+\phi(t,x,y,z)u_{xz},\nonumber\\
W&=&-u_{x}.
\end{eqnarray}
For $\Gamma_{42}$ the components of the conserved vectors are
\begin{eqnarray}
C^{t}&=&D_{t}^{1-\alpha}W\phi(t,x,y,z)+J(W,\phi_t),\nonumber\\
C^{x}&=&W\phi_x+\phi(t,x,y,z)u_{xy},\nonumber\\
C^{y}&=&\phi(t,x,y,z)(D_{t}^{1-\alpha}u-(u_{xx}+u_{yy}+u_{zz}))+W\phi_y+\phi(t,x,y,z)u_{yy},\nonumber\\
C^{z}&=&W\phi_z+\phi(t,x,y,z)u_{yz},\nonumber\\
W&=&-u_{y}.
\end{eqnarray}
For $\Gamma_{43}$ the components of the conserved vectors are
\begin{eqnarray}
C^{t}&=&D_{t}^{1-\alpha}(W)\phi(t,x,y,z)+J(W,\phi_t),\nonumber\\
C^{x}&=&W\phi_x+\phi(t,x,y,z)u_{xz},\nonumber\\
C^{y}&=&W\phi_y+\phi(t,x,y,z)u_{yz},\nonumber\\
C^{z}&=&\phi(t,x,y,z)(D_{t}^{1-\alpha}u-(u_{xx}+u_{yy}+u_{zz}))+W\phi_z+\phi(t,x,y,z)u_{zz},\nonumber\\
W&=&-u_{z}.
\end{eqnarray}
For $\Gamma_{44}$ the components of the conserved vectors are
\begin{eqnarray}
C^{t}&=&D_{t}^{1-\alpha}W\phi(t,x,y,z)+J(W,\phi_t),\nonumber\\
C^{x}&=&-y\phi(t,x,y,z)(D_{t}^{1-\alpha}u-(u_{xx}+u_{yy}+u_{zz}))+W\phi_x-\phi(t,x,y,z)(yu_{xx}-(xu_{xy}+u_{y})),\nonumber\\
C^{y}&=&x\phi(t,x,y,z)(D_{t}^{1-\alpha}u-(u_{xx}+u_{yy}+u_{zz}))+W\phi_y-\phi(t,x,y,z)((u_{x}+yu_{xy})-xu_{yy}),\nonumber\\
C^{z}&=&W\phi_z+\phi(t,x,y,z)(yu_{xz}-xu_{yz}),\nonumber\\
W&=&yu_{x}-xu_{y}.\nonumber\\
\end{eqnarray}
For $\Gamma_{45}$ the components of the conserved vectors are
\begin{eqnarray}
C^{t}&=&D_{t}^{1-\alpha}(W)\phi(t,x,y,z)+J(W,\phi_t),\nonumber\\
C^{x}&=&W\phi_x+\phi(t,x,y,z)(yu_{xz}-zu_{xy}),\nonumber\\
C^{y}&=&z\phi(t,x,y,z)(D_{t}^{1-\alpha}u-(u_{xx}+u_{yy}+u_{zz}))+W\phi_y-\phi(t,x,y,z)((u_{z}+yu_{zy})-zu_{yy}),\nonumber\\
C^{z}&=&-y\phi(t,x,y,z)(D_{t}^{1-\alpha}u-(u_{xx}+u_{yy}+u_{zz}))+W\phi_z-\phi(t,x,y,z)(yu_{zz}-(u_{y}+zu_{yz})),\nonumber\\
W&=&yu_{z}-zu_{y}.
\end{eqnarray}
For  $\Gamma_{46}$ the components of the conserved vectors are
\begin{eqnarray}
C^{t}&=&D_{t}^{1-\alpha}(W)\phi(t,x,y,z)+J(W,\phi_t),\nonumber\\
C^{x}&=&z\phi(t,x,y,z)(D_{t}^{1-\alpha}u-(u_{xx}+u_{yy}+u_{zz}))+W\phi_x-\phi(t,x,y,z)((u_{z}+xu_{xz})-zu_{xx}),\nonumber\\
C^{y}&=&W\phi_{y}+\phi(t,x,y,z)(xu_{yz}-zu_{xy}),\nonumber\\
C^{z}&=&-x\phi(t,x,y,z)(D_{t}^{1-\alpha}u-(u_{xx}+u_{yy}+u_{zz}))+W\phi_z-\phi(t,x,y,z)(xu_{zz}-(u_{x}+u_{xz})),\nonumber\\
W&=&xu_{z}-zu_{x}.
\end{eqnarray}
For $\Gamma_{47}$ the components of the conserved vectors are
\begin{align}
C^{t}&=\begin{aligned}[t]
&2t\phi(t,x,y,z)(D_{t}^{1-\alpha}u-(u_{xx}+u_{yy}+u_{zz}))+D_{t}^{1-\alpha}(W)\phi(t,x,y,z)+J(W,\phi_t),\nonumber\\
\end{aligned}\\
C^{x}&=\begin{aligned}[t]
&\alpha x\phi(t,x,y,z)(D_{t}^{1-\alpha}u-(u_{xx}+u_{yy}+u_{zz}))+W\phi_{x},\nonumber\\
&+\phi(t,x,y,z)(u_x(\alpha-1)-2tu_{xt}-\alpha(u_x+xu_{xx})-\alpha yu_{xy}-\alpha zu_{xz}),\nonumber\\
\end{aligned}\\
C^{y}&=\begin{aligned}[t]
&\alpha y\phi(t,x,y,z)(D_{t}^{1-\alpha}u-(u_{xx}+u_{yy}+u_{zz}))+W\phi_{y}+\nonumber\\
&+\phi(t,x,y,z)(u_{y}(\alpha-1)-2tu_{yt}-\alpha(u_{y}+yu_{yy})-\alpha xu_{xy}-\alpha zu_{yz})\nonumber\\
\end{aligned}\\
C^{z}&=\begin{aligned}[t]
&\alpha z\phi(t,x,y,z)(D_{t}^{1-\alpha}u-(u_{xx}+u_{yy}+u_{zz}))+W\phi_{z}+\nonumber\\
&+\phi(t,x,y,z)(u_{z}(\alpha-1)-2tu_{zt}-\alpha(u_{z}+zu_{zz})-\alpha xu_{xz}-\alpha yu_{yz})\nonumber\\
\end{aligned}\\
W&=\begin{aligned}[t]
&u(\alpha-1)-2tu_{t}-\alpha x u_{x}-\alpha y u_{y}-\alpha z u_{z}.\nonumber\\
\end{aligned}\\
\end{align}
For $\Gamma_{48}$ the components of the conserved vectors are
\begin{eqnarray}
C^{t}&=&D_{t}^{1-\alpha}(W)\phi(t,x,y,z)+J(W,\phi_t),\nonumber\\
C^{x}&=&W\phi_x-\phi(t,x,y,z)u_{x},\nonumber\\
C^{y}&=&W\phi_y-\phi(t,x,y,z)u_{y},\nonumber\\
C^{z}&=&W\phi_z-\phi(t,x,y,z)u_{z},\nonumber\\
W&=&u.
\end{eqnarray}
For $\Gamma_{49}$ the components of the conserved vectors are
\begin{eqnarray}
C^{t}&=&D_{t}^{1-\alpha}(W)\phi(t,x,y,z)+J(W,\phi_t),\nonumber\\
C^{x}&=&W\phi_x-\phi(t,x,y,z)F_x,\nonumber\\
C^{y}&=&W\phi_y-\phi(t,x,y,z)F_y,\nonumber\\
C^{z}&=&W\phi_z-\phi(t,x,y,z)F_z,\nonumber\\
W&=&F(t,x,y,z).
\end{eqnarray}

\subsection{The conservation laws for the four-dimensional time-fractional heat equation}

The components of the conserved vectors for $\Gamma_{61}$ are
\begin{eqnarray}
C^{t}&=&D_{t}^{1-\alpha}(W)\phi(t,x,y,z,w)+J(W,\phi_t),\nonumber\\
C^{x}&=&\phi(t,x,y,z,w)(D_{t}^{1-\alpha}u-(u_{xx}+u_{yy}+u_{zz}+u_{ww}))+W\phi_x+\phi(t,x,y,z,w) u_{xx},\nonumber\\
C^{y}&=&W\phi_{y}+\phi(t,x,y,z,w)u_{xy},\nonumber\\
C^{z}&=&W\phi_{z}+\phi(t,x,y,z,w)u_{xz},\nonumber\\
C^{w}&=&W\phi_{w}+\phi(t,x,y,z,w)u_{xw},\nonumber\\
W&=&-u_x.\nonumber\\
\end{eqnarray}
The components of the conserved vectors for $\Gamma_{62}$ are
\begin{eqnarray}
C^{t}&=&D_{t}^{1-\alpha}(W)\phi(t,x,y,z,w)+J(W,\phi_t),\nonumber\\
C^{x}&=&W\phi_{x}+\phi(t,x,y,z,w)u_{xy},\nonumber\\
C^{y}&=&\phi(t,x,y,z,w)(D_{t}^{1-\alpha}u-(u_{xx}+u_{yy}+u_{zz}+u_{ww}))+W\phi_y+\phi(t,x,y,z,w) u_{yy},\nonumber\\
C^{z}&=&W\phi_{z}+\phi(t,x,y,z,w)u_{zy},\nonumber\\
C^{w}&=&W\phi_{w}+\phi(t,x,y,z,w)u_{yw},\nonumber\\
W&=&-u_y.
\end{eqnarray}
The components of the conserved vectors for $\Gamma_{63}$ are
\begin{eqnarray}
C^{t}&=&D_{t}^{1-\alpha}(W)\phi(t,x,y,z,w)+J(W,\phi_t),\nonumber\\
C^{x}&=&W\phi_{x}+\phi(t,x,y,z,w)u_{xz},\nonumber\\
C^{y}&=&W\phi_{y}+\phi(t,x,y,z,w)u_{yz},\nonumber\\
C^{z}&=&\phi(t,x,y,z,w)(D_{t}^{1-\alpha}u-(u_{xx}+u_{yy}+u_{zz}+u_{ww}))+W\phi_z+\phi(t,x,y,z,w) u_{zz},\nonumber\\
C^{w}&=&W\phi_{w}+\phi(t,x,y,z,w)u_{wz},\nonumber\\
W&=&-u_z.\nonumber\\
\end{eqnarray}
The components of the conserved vectors for $\Gamma_{64}$ are
\begin{eqnarray}
C^{t}&=&D_{t}^{1-\alpha}(W)\phi(t,x,y,z,w)+J(W,\phi_t),\nonumber\\
C^{x}&=&W\phi_{x}+\phi(t,x,y,z,w)u_{xw},\nonumber\\
C^{y}&=&W\phi_{y}+\phi(t,x,y,z,w)u_{yw},\nonumber\\
C^{z}&=&W\phi_{z}+\phi(t,x,y,z,w)u_{zw},\nonumber\\
C^{w}&=&\phi(t,x,y,z,w)(D_{t}^{1-\alpha}u-(u_{xx}+u_{yy}+u_{zz}+u_{ww}))+W\phi_w+\phi(t,x,y,z,w) u_{ww},\nonumber\\
W&=&-u_w.
\end{eqnarray}
The components of the conserved vectors for $\Gamma_{65}$ are
\begin{eqnarray}
C^{t}&=&D_{t}^{1-\alpha}(W)\phi(t,x,y,z,w)+J(W,\phi_t),\nonumber\\
C^{x}&=&-y\phi(t,x,y,z,w)(D_{t}^{1-\alpha}u-(u_{xx}+u_{yy}+u_{zz}+u_{ww}))+W\phi_{x}-\phi(yu_{xx}-(u_{xy}+u_{y})),\nonumber\\
C^{y}&=&x\phi(t,x,y,z,w)(D_{t}^{1-\alpha}u-(u_{xx}+u_{yy}+u_{zz}+u_{ww}))+W\phi_y-\phi((u_{x}+u_{xy})-xu_{yy}),\nonumber\\
C^{z}&=&W\phi_{z}+\phi(t,x,y,z,w)(yu_{xz}-xu_{yz}),\nonumber\\
C^{w}&=&W\phi_{w}+\phi(t,x,y,z,w)(yu_{xw}-xu_{yw}),\nonumber\\
W&=&yu_{x}-xu_{y}.
\end{eqnarray}
The components of the conserved vectors for $\Gamma_{66}$ are
\begin{eqnarray}
C^{t}&=&D_{t}^{1-\alpha}(W)\phi(t,x,y,z,w)+J(W,\phi_t),\nonumber\\
C^{x}&=&W\phi_x+\phi(t,x,y,z,w)(yu_{xz}-zu_{xy}),\nonumber\\
C^{y}&=&z\phi(t,x,y,z,w)(D_{t}^{1-\alpha}u-(u_{xx}+u_{yy}+u_{zz}+u_{ww}))+W\phi_{y}-\phi(t,x,y,z,w)((u_{z}+u_{zy})-zu_{yy}),\nonumber\\
C^{z}&=&-y\phi(t,x,y,z,w)(D_{t}^{1-\alpha}u-(u_{xx}+u_{yy}+u_{zz}+u_{ww}))+W\phi_z-\phi(t,x,y,z,w)(yu_{zz}-(u_{y}+u_{yz})),\nonumber\\
C^{w}&=&W\phi_{w}+\phi(t,x,y,z,w)(yu_{zw}-zu_{yw}),\nonumber\\
W&=&yu_{z}-zu_{y}.
\end{eqnarray}
The components of the conserved vectors for $\Gamma_{67}$ are
\begin{eqnarray}
C^{t}&=&D_{t}^{1-\alpha}(W)\phi(t,x,y,z,w)+J(W,\phi_t),\nonumber\\
C^{x}&=&W\phi_{x}+\phi(t,x,y,z,w)(wu_{xy}-yu_{xw}),\nonumber\\
C^{y}&=&-w\phi(t,x,y,z,w)(D_{t}^{1-\alpha}u-(u_{xx}+u_{yy}+u_{zz}+u_{ww}))+W\phi_y-\phi(t,x,y,z,w)(wu_{yy}-(yu_{wy}+u_w)),\nonumber\\
C^{z}&=&W\phi_{z}+\phi(t,x,y,z,w)(wu_{yz}-yu_{wz}),\nonumber\\
C^{w}&=&y\phi(t,x,y,z,w)(D_{t}^{1-\alpha}u-(u_{xx}+u_{yy}+u_{zz}+u_{ww}))+W\phi_w-\phi(t,x,y,z,w)((wu_{yw}+u_{y})-yu_{ww}),\nonumber\\
W&=&wu_y-yu_w.
\end{eqnarray}
The components of the conserved vectors for $\Gamma_{69}$ are
\begin{eqnarray}
C^{t}&=&D_{t}^{1-\alpha}(W)\phi(t,x,y,z,w)+J(W,\phi_t),\nonumber\\
C^{x}&=&z\phi(t,x,y,z,w)(D_{t}^{1-\alpha}u-(u_{xx}+u_{yy}+u_{zz}+u_{ww}))+W(\phi_x)-\phi(t,x,y,z,w)((xu_{xz}+u_{z})-zu_{xx}),\nonumber\\
C^{y}&=&W\phi_{y}+\phi(t,x,y,z,w)(xu_{zy}-zu_{xy}),\nonumber\\
C^{z}&=&-x\phi(t,x,y,z,w)(D_{t}^{1-\alpha}u-(u_{xx}+u_{yy}+u_{zz}+u_{ww}))+W\phi_z-\phi(t,x,y,z,w)(xu_{zz}-(zu_{xz}+u_x)),\nonumber\\
C^{w}&=&W\phi_{w}+\phi(t,x,y,z,w)(xu_{zw}-zu_{xw}),\nonumber\\
W&=&xu_{z}-zu_{x}.
\end{eqnarray}
The components of the conserved vectors for $\Gamma_{610}$ are
\begin{eqnarray}
C^{t}&=&D_{t}^{1-\alpha}(W)\phi(t,x,y,z,w)+J(W,\phi_t),\nonumber\\
C^{x}&=&-w\phi(t,x,y,z,w)(D_{t}^{1-\alpha}u-(u_{xx}+u_{yy}+u_{zz}+u_{ww}))+W\phi_x-\phi(t,x,y,z,w)(wu_{xx}-(xu_{xw}+u_w)),\nonumber\\
C^{y}&=&W\phi_{y}+\phi(t,x,y,z,w)(wu_{xy}-xu_{wy}),\nonumber\\
C^{z}&=&W\phi_{z}+\phi(t,x,y,z,w)(wu_{xz}-xu_{wz}),\nonumber\\
C^{w}&=&x\phi(t,x,y,z,w)(D_{t}^{1-\alpha}u-(u_{xx}+u_{yy}+u_{zz}+u_{ww}))+W\phi_w-\phi((wu_{xw}+u_{x})-xu_{ww}),\nonumber\\
W&=&wu_x-xu_w.
\end{eqnarray}
The components of the conserved  vectors for $\Gamma_{611}$ are
\begin{eqnarray}
C^{t}&=&D_{t}^{1-\alpha}(W)\phi(t,x,y,z,w)+J(W,\phi_t),\nonumber\\
C^{x}&=&W\phi_{x}+\phi(t,x,y,z,w)(wu_{xz}-zu_{xw}),\nonumber\\
C^{y}&=&W\phi_{y}+\phi(t,x,y,z,w)(wu_{yz}-zu_{yw}),\nonumber\\
C^{z}&=&-w\phi(t,x,y,z,w)(D_{t}^{1-\alpha}u-(u_{xx}+u_{yy}+u_{zz}+u_{ww}))+W\phi_z-\phi(t,x,y,z,w)(wu_{zz}-(zu_{wz}+u_w)),\nonumber\\
C^{w}&=&z\phi(t,x,y,z,w)(D_{t}^{1-\alpha}u-(u_{xx}+u_{yy}+u_{zz}+u_{ww}))+W\phi_w-\phi((u_z+u_{zw})-zu_{ww}),\nonumber\\
W&=&wu_{z}-zu_{w}.
\end{eqnarray}
The components of the conserved vectors for $\Gamma_{612}$ are
\begin{align}
C^{t}&=\begin{aligned}[t]
&2t\phi(t,x,y,z,w)(D_{t}^{1-\alpha}u-(u_{xx}+u_{yy}+u_{zz}+u_{ww}))+D_{t}^{1-\alpha}(W)\phi(t,x,y,z)+J(W,\phi_t),\nonumber\\
\end{aligned}\\
C^{x}&=\begin{aligned}[t]
&\alpha x\phi(t,x,y,z,w)(D_{t}^{1-\alpha}u-(u_{xx}+u_{yy}+u_{zz}+u_{ww}))+W\phi_x-\phi(t,x,y,z,w)((\alpha-1)u_x\nonumber\\
&-\alpha(u_x+xu_{xx})-\alpha yu_{xy}-\alpha zu_{xz}-\alpha wu_{xw}-2tu_{xt}),\nonumber\\
\end{aligned}\\
C^{y}&=\begin{aligned}[t]
&\alpha y\phi(t,x,y,z,w)(D_{t}^{1-\alpha}u-(u_{xx}+u_{yy}+u_{zz}+u_{ww}))+W\phi_y-\phi(t,x,y,z,w)((\alpha-1)u_y\nonumber\\
&-\alpha(u_y+yu_{yy})-\alpha xu_{xy}-\alpha zu_{zy}-\alpha wu_{wy}-2tu_{yt}),\nonumber\\
\end{aligned}\\
C^{z}&=\begin{aligned}[t]
&\alpha z\phi(t,x,y,z,w)(D_{t}^{1-\alpha}u-(u_{xx}+u_{yy}+u_{zz}+u_{ww}))+W\phi_z-\phi(t,x,y,z,w)((\alpha-1)u_z\nonumber\\
&-\alpha xu_{xz}-\alpha yu_{yz}-\alpha(u_z+zu_{zz})-\alpha wu_{wz}-2tu_{zt}),\nonumber\\
\end{aligned}\\
C^{w}&=\begin{aligned}[t]
&\alpha w\phi(t,x,y,z,w)(D_{t}^{1-\alpha}u-(u_{xx}+u_{yy}+u_{zz}+u_{ww}))+W\phi_w-\phi(t,x,y,z,w)((\alpha-1)u_w\nonumber\\
&-\alpha xu_{xw}-\alpha yu_{yw}-\alpha zu_{zw}-\alpha(u_w+wu_{ww})-2tu_{wt}),\nonumber\\
\end{aligned}\\
W&=\begin{aligned}[t]
&u(\alpha-1)-\alpha x u_{x}-\alpha y u_{y}-\alpha z u_{z}-\alpha w u_{w}-2t u_{t}.\nonumber\\
\end{aligned}\\
\end{align}
For $\Gamma_{613}$ the components of the conserved vectors are
\begin{eqnarray}
C^{t}&=&D_{t}^{1-\alpha}(W)\phi(t,x,y,z,w)+J(W,\phi_t),\nonumber\\
C^{x}&=&W\phi_x-\phi(t,x,y,z,w)u_{x},\nonumber\\
C^{y}&=&W\phi_y-\phi(t,x,y,z,w)u_{y},\nonumber\\
C^{z}&=&W\phi_z-\phi(t,x,y,z,w)u_{z},\nonumber\\
C^{w}&=&W\phi_w-\phi(t,x,y,z,w)u_{w},\nonumber\\
W&=&u.
\end{eqnarray}
For $\Gamma_{614}$ the components of the conserved vectors are
\begin{eqnarray}
C^{t}&=&D_{t}^{1-\alpha}(W)\phi(t,x,y,z,w)+J(W,\phi_t),\nonumber\\
C^{x}&=&W\phi_x-\phi(t,x,y,z,w)F_x,\nonumber\\
C^{y}&=&W\phi_y-\phi(t,x,y,z,w)F_y,\nonumber\\
C^{z}&=&W\phi_z-\phi(t,x,y,z,w)F_z,\nonumber\\
C^{w}&=&W\phi_w-\phi(t,x,y,z,w)F_w,\nonumber\\
W&=&F(t,x,y,z,w).
\end{eqnarray}

\subsection{The conservation laws for the n-dimensional time-fractional heat equation}

The conservations laws for each of the symmetry $\Gamma_{7j}$, $j=1,2,..5$  can be summarised as:
\begin{eqnarray}
C^{t}&=&\xi^{0}L+D_{t}^{1-\alpha}(W)\phi(t,x,y,z,w)+J(W,\phi_t),\nonumber\\
C^{x_{i}}&=&\xi^{i}L+W\left(\frac{\partial L}{\partial_{u_{x_{i}}}}-D_{x_{i}}\frac{\partial L}{\partial_{u_{x_{i}x_{i}}}}\right)+D_{x_{i}}(W)\frac{\partial L}{\partial_{u_{x_{i}x_{i}}}}, i=1,2,....n.\nonumber\\
\end{eqnarray}
The terms $L$, $W$, $\xi^{0}$ and $\xi^{i}$ for various integral values of $i$, are already defined in the preliminaries.

\section{Concluding remarks and discussion}

We have investigated time-fractional heat equation using Lie symmetries and obtaining a classification of these symmetries.

We note that for all dimensions, when we consider the case $0<\alpha<1$, we lose the translational symmetry, $\partial_{t}$, as explained in
Gazizov {\it et al} (2009).

The number of symmetries is reduced significantly and the differences in Lie Algebras for the fractional- and integeral-order PDEs can be attributed to this fact.
The significance of the reduction in conservation laws for the fractional form of the nonfractional case can be a matter of further research..

According to Myeni and Leach (2009) in the case of linear ODEs the number of solution symmetries is equal to the order of the equation.
From this paper we see that for integer-order linear PDEs the number of solution symmetries is equal to the product of the order and space dimension,
whereas for the fractional PDEs it is half of the product of the order and space dimension.

We have generalised the number of symmetries we can find for an $n$-dimensional time-fractional heat equation. For the case of $\alpha=1$
the number of symmetries for the $n$-dimensional case is, $$\frac{1}{2}(n^{2}+3n+10).$$ In the case of $0<\alpha<1$ the number of symmetries is,
$$\frac{1}{2}(n^{2}+n+6).$$

\section*{Acknowledgements}

 AKH expresses grateful thanks to UGC (India), NFSC, Award No. F1-17.1/201718/RGNF-2017-18-SC-ORI-39488 for financial support. PGLL thanks the University of
Kwazulu-Natal, Durban University of
Technology, and the National Research Foundation of South Africa for financial support and to the Department of Mathematics, Pondicherry University, for
gracious hospitality.

\end{document}